\documentclass[reqno]{amsart}
\usepackage{amssymb}
\usepackage[usenames, dvipsnames]{color}
\usepackage{txfonts}
\usepackage{hyperref}
\usepackage{verbatim}
\usepackage[pagewise]{lineno}
\usepackage{marginnote}
\usepackage{todonotes}
\usepackage[normalem]{ulem}
\usepackage{cancel}
\usepackage{soul}

\theoremstyle{plain}
\newtheorem{theorem}{Theorem}[section]
\newtheorem{lemma}[theorem]{Lemma}
\newtheorem{corollary}[theorem]{Corollary}
\newtheorem{proposition}[theorem]{Proposition}

\theoremstyle{definition}
\newtheorem{definition}[theorem]{Definition}

\theoremstyle{remark}
\newtheorem{remark}[theorem]{Remark}

\newcommand{\pf}{\noindent {\bf Proof. \hspace{2mm}}}

\def\p{\partial}

\def\be{\begin{equation}}
\def\ee{\end{equation}}
\def\bes{\begin{equation*}}
\def\ees{\end{equation*}}
\def\bali{\begin{aligned}}
	\def\eali{\end{aligned}}

\numberwithin{equation}{section}

\makeatletter
\def\dashint{\operatorname%
	{\,\,\text{\bf--}\kern-.98em\DOTSI\intop\ilimits@\!\!}}
\makeatother

\begin{document}
	
	\title[Time analyticity]{Time analyticity of the biharmonic heat equation, the heat equation with potentials and some nonlinear heat equations}

	\author[C. Zeng]{Chulan Zeng}
	\address[C. Zeng]{Department of Mathematics, University of California, Riverside, CA, 92521, USA}
	
	\email{czeng011@ucr.edu}
	
	\thanks{}
	
	\subjclass[2010]{}
	

	\date{\today}
	
	
	\begin{abstract}
		In this paper, we investigate the pointwise time analyticity of three differential equations. They are the biharmonic heat equation, the heat equation with potentials and some nonlinear heat equations with power nonlinearity of order $p$. The potentials include all the nonnegative ones. For the first two equations, we prove if $u$ satisfies some growth conditions in $(x,t)\in \mathrm{M}\times [0,1]$, then $u$ is analytic in time $(0,1]$. Here $\mathrm{M}$ is $R^d$
		or a complete noncompact manifold with Ricci curvature bounded from below by a
		constant. Then we obtain a necessary and sufficient condition such that $u(x,t)$ is analytic in time at $t=0$. Applying this method, we also obtain a necessary and sufficient condition for the solvability of the backward equations, which is ill-posed in general.
		\par  For the nonlinear heat equation with power nonlinearity of order $p$, we prove that a solution is analytic in time $t\in (0,1]$ if it is bounded in $\mathrm{M}\times[0,1]$ and $p$ is a positive integer. In addition, we investigate the case when $p$ is a rational number with a stronger assumption $0<C_3 \leq |u(x,t)| \leq C_4$. It is also shown that a solution may not be analytic in time if it is allowed to be $0$. As a lemma, we obtain an estimate of $\p_t^k \Gamma(x,t;y)$ where $\Gamma(x,t;y)$ is the heat kernel on a manifold, with an explicit estimation of the coefficients.
		\par An interesting point is that a solution may be analytic in time even if it is not smooth in the space variable $x$, implying that the analyticity of space and time can be independent. Besides, for general manifolds, space analyticity may not hold since it requires certain bounds on curvature
		and its derivatives.\\
		
	\end{abstract}

	\maketitle
	
	\tableofcontents

	\section{Introduction}

	In this paper, we investigate the pointwise time analyticity of three differential equations. The first one is the biharmonic heat equation  
	\be\label{2ndeq}
	\partial_{t} u+\Delta^2 u=0, \quad\forall (x,t)\in \mathrm{M}\times [0,1].
	\ee 
	
	Here and below, $\mathrm{M}$ is $R^d$
	or a $d$ demensional complete noncompact manifold with Ricci curvature bounded from below by a
	constant.
	The second one is the heat equation with potentials
	\be\label{3rdeq} 
	\p_tu(x,t)-\Delta u(x,t)+V(x)u(x,t)=0, \quad\forall (x,t)\in \mathrm{M}\times [0,1],\ d\geq 3.
	\ee
	In one case, $V=V(x)$ is a potential function in $L^q(\mathrm{M})$ for some $q \geq 1$, with some growth conditions. In another case, we treat $V(x)\geq0$. The last equation is some nonlinear heat equations with power nonlinearity of order $p$ where $p$ is some positive rational number,
	\begin{equation}\label{equa0}
	u_t(x,t)-\Delta u(x,t)=u^p(x,t),\quad\forall (x,t)\in \mathrm{M}\times [0,1].
	\end{equation}
	The goal of this paper is to extend the results in H.Dong$\&$Q.Zhang\cite{[Zhang]} to these three differential equations above. 
	
	While the spatial analyticity is usually
	true for generic solutions, the time analyticity is harder to prove and is false in general.
	For example, 
	it is not difficult 
	to construct a solution of
	the heat equation in a 
	space-time cylinder in the Euclidean setting, which is not time analytic in a sequence of
	moments. 
	Besides, the time analyticity is not a local property, so we need to impose certain growth
	conditions on solutions and data at infinity. Under various assumptions, there
	are numerous time-analyticity results for the heat equation and other parabolic type
	equations. See, for example, D.Widder\cite{[Widder]} and H.Dong$\&$Q.Zhang\cite{[Zhang]}. Moreover, if one
	imposes zero boundary conditions on the lateral boundary of a smooth cylindrical domain,
	then certain solutions of the heat, biharmonic heat, and many other parabolic equations are
	analytic in time. See, for example, K. Masuda\cite{[Masuda]}, G. Komatsu\cite{[Komatsu]}, Y.Giga\cite{[GI]}, and L.Escauriaza, S.Montaner$\&$C.Zhang\cite{[LE]}. One can also consider
	solutions in certain $L_p$ spaces with $p\in(1,\infty)$, see K.Promislow\cite{[Pr]} for large class of dissipative equations in the periodic setting.
	
	In a related development, there have been increasing interest in the study of time analyticity of solutions of parabolic type
	equations on the Euclidean and manifold setting. For example, in the papers H.Dong$\&$Q.Zhang\cite{[Zhang]}, it is proven that if a smooth solution of the heat equation in $\mathrm{M}\times (-2,0]$ is of exponential growth of order 2, then it is analytic in time in $t\in[-1,0]$. In Q.Zhang\cite{[Zhang2]}, it is pointed out that the time analyticity is equivalent to time inversibility.
	Besides,
	it is proven in L.Escauriaza, S.Montaner and C.Zhang\cite{[LE]} that for any bounded domain $\Omega\subset R^d$ with analytic boundary, 
	any solution of the high order harmonic heat equation
	$$
	\quad\left\{\begin{array}{l}
	u_t+(-\Delta)^m u=0, \quad\forall (x,t)\in\Omega\times(0,1], \\
	u=Du=\cdots=D^{m-1}u=0\ \text{on}\ \partial\Omega\times (0,1],\ u(x,0)\in L^2(\Omega)
	\end{array}\right.
	$$
	is analytic in time $t\in(0,1]$.
	There are also some other results about time analyticity of parabolic type differential equations with noncompact boundary conditions in H.Dong$\&$X.Pan\\\cite{[DP]}.
	
	Here are the main results of this paper. The first one is about the biharmonic heat equation \eqref{2ndeq}.
	\begin{theorem}\label{2ndtheo}
		Let $\mathrm{M}$ be a d dimensional, complete, noncompact Riemannian manifold such that the Ricci curvature satisfies $Ric \geq-(d - 1)K_0$ for a nonnegative constant $K_0$. 
		
		Let $u=u(x,t)$ be a smooth solution of the biharmonic heat equation \eqref{2ndeq} on $\mathrm{M}\times [0,1]$ of exponential growth of order $\frac{4}{3},$ namely
		
		$$|u(x,t)| \leq A_{1} e^{A_{2} d^{\frac{4}{3}}(x,0)}, \quad \forall(x,t) \in\mathrm{M}\times [0,1] , $$
		where $A_{1}$ and $A_{2}$ are positive constants. Then $u$ is analytic in time $t \in(0,1]$ with radius of convergence depending only on $t$, $d$, $K_0$ and $A_2$. Moreover, if $t\in(1-\delta,1]$ for some small $\delta>0$, we have
		$$
		u(x,t)=\sum_{j=0}^{\infty} a_{j}(x) \frac{(t-1)^{j}}{j !}
		$$
		with $-\Delta^2 a_{j}(x)=a_{j+1}(x),$ and
		$$
		\left|a_{j}(x)\right|=\left|(-\Delta^2)^j a_{0}(x)\right| \leq  A^* A_{3}^{j+1} j^{j} e^{2 A_{2} d^{\frac{4}{3}}(x, 0)}, \quad j=0,1,2, \ldots
		$$
		where $A_{3}=A_3(d,K_0,A_2)$ and $A^*=A^*(A_1,d,x_0,\mathrm{M})$.
	\end{theorem}
	Then we have two main theorems about the heat equation with potentials \eqref{3rdeq}. We define the weak solution in the beginning of Section 3.
	\begin{theorem}\label{3rdtheo}
		Let $\mathrm{M}$ be a $d$ dimensional, complete, noncompact, smooth Riemannian manifold such that the Ricci curvature satisfies $Ric \geq-(d - 1)K_0$ for some nonnegative constant $K_0$ and 
		$$\inf\limits_{x\in\mathrm{M}}|B(x,1)|>0.$$
		Assume $V=V(x)$ satisfies the following conditions:\\
		(1) There exists some $R^*>0$ such that $V(\cdot)\in L^{q}(B(0,R^*))$ for some $q>\frac{d}{2}$. \\
		(2) For some constant $C^{**}>0$, if $d(x,0)>R^*$, then $|V(x)|\leq C^{**} d(x,0)^{\alpha}$ where $\alpha=\frac{2q-d}{q-1}$ and $d> 2$.\\
		(3) $V(\cdot)\in L^1(\mathrm{M}\backslash B(0,R^*))$ and assume $\|V\|_{L^1(\mathrm{M}\backslash B(0,R^*))}=D^*$.
		
		Let $$\|V\|_{L^{q}(B(0,R^*))}=C^{*}$$ where $C^*$ is a positive constant and let $u=u(x,t)$ be a weak solution of equation \eqref{3rdeq} on $\mathrm{M}\times [0,1]$ of exponential growth of order $2,$ namely
		$$|u(x,t)| \leq A_{1} e^{A_{2} d^{2}(x, 0)}, \quad \forall(x,t) \in\mathrm{M}\times [0,1] , $$
		where $A_{1}$ and $A_{2}$ are some positive constants. Then $u$ is analytic in $t \in (0,1/2]$ with radius of convergence depending only on $t$, $d$, $q$, $K_0$, $A_{2}$, $\alpha$ and $C^*$.
		
		Moreover, if $t\in(1/2-\delta,1/2]$ for some small $\delta>0$, we have
		$$
		u(x,t)=\sum_{j=0}^{\infty} a_{j}(x) \frac{(t-1/2)^{j}}{j !}
		$$
		with $(\Delta-V) a_{j}(x)=a_{j+1}(x),$ and
		\be\label{condnew}
		\left|a_{j}(x)\right| =\left|(\Delta-V)^j a_0(x)\right|\leq A_1A_{3}^{j+1} j^{j} e^{A_4 d^{2}(x, 0)}, \quad j=0,1,2, \ldots
		\ee
		where constants $A_3=A_3(d,q,K_0,A_2,\alpha,C^*)$ and $A_4=A_4(A_2,\alpha,C^{**},D^*)$.\\
	\end{theorem}
	Here the extra condition $d\geq 3$ can be removed in the case of $R^d$. We will explain in more detail during the proof.
	\begin{theorem}\label{4ththeo}
		Let $\mathrm{M}$ be a $d$ dimensional, complete, noncompact Riemannian manifold such that the Ricci curvature satisfies $Ric \geq-(d - 1)K_0$ for some nonnegative constant $K_0$. 
		
		Let $u=u(x,t)$ be a weak solution of the heat equation with nonnegative potentials \eqref{3rdeq} where $V=V(x)\geq 0$ on $\mathrm{M}\times [0,1]$. If $u$ is of exponential growth of order $2$, namely
		$$|u(x,t)| \leq A_{1} e^{A_{2} d^{2}(x, 0)}, \quad \forall(x,t) \in\mathrm{M}\times [0,1] , $$
		where $A_{1}$ and $A_{2}$ are positive constants, then $u$ is analytic in $t \in(0,1]$ with radius depending only on $t$, $d$, $K_0$ and $A_{2}$. 
		
		Moreover, if $t\in(1-\delta,1]$ for some small $\delta>0$, we have
		$$
		u(x,t)=\sum_{j=0}^{\infty} a_{j}(x) \frac{(t-1)^{j}}{j !}
		$$
		with $(\Delta-V) a_{j}(x)=a_{j+1}(x),$ and
		\be\label{connew2}
		\left|a_{j}(x)\right| =\big|\left(\Delta-V \right)^j a_0(x)\big|\leq  A_1A_{5}^{j+1} j^{j} e^{2 A_{2} d^{2}(x, 0)}, \quad j=0,1,2, \ldots
		\ee
		where $A_{5}=A_{5}(d, K_0,A_{2})$.
	\end{theorem}
	We also have two theorems about some nonlinear heat equations with power nonlinearity of order $p$.
	\begin{theorem}\label{theo00}
		Let $\mathrm{M}$ be a d dimensional, complete, noncompact Riemannian manifold such that the Ricci curvature satisfies $Ric \geq-(d - 1)K_0$ for some nonnegative constant $K_0$. 
		
		Let $u=u(x,t)$ be a solution to equation \eqref{equa0} where $p$ is a positive integer. Suppose $u$
		satisfies 
		$$|u(x,t)|\leq C_2 \quad\text{in}\quad \mathrm{M}\times [0,1],$$ 
		for some constant $C_2$.
		Then $u$ is analytic in time for any $t\in(0,1]$ with radius of convergence independent of $x$.\\
	\end{theorem}
	\begin{theorem}\label{theo05}
		Let $\mathrm{M}$ be the same manifold as Theorem \ref{theo00} above and $p=q_1/q_2$ for some positive integers $q_1$ and $q_2$.
		Assume that a solution $u=u(x,t)$ to the equation \eqref{equa0} satisfies 
		$$0<C_3 \leq |u(x,t)| \leq C_4 \quad\text{in}\quad \mathrm{M}\times [0,1],$$ 
		where $C_3$, $C_4$ are some constants. Then $u$ is analytic in time for any $t\in(0,1]$ with radius of convergence independent of $x$. 
	\end{theorem}
	
	Now we give a brief outline of this paper. In Section \ref{2ndsec}, we investigate the time analyticity of the biharmonic heat equation \eqref{2ndeq}. As a corollary, we obtain a necessary and sufficient condition for the solvability of the backward biharmonic heat equation $\partial_{t} u-\Delta^2 u=0$. As another corollary, we also obtain a necessary and sufficient condition under which the solution of \eqref{2ndeq} is analytic in time at initial time $t=0$. Section \ref{3rdsec} pertains the time analyticity of the heat equation with potentials \eqref{3rdeq}. We use similar methods and obtain similar results as in Section \ref{2ndsec}. We demonstrate some solutions which may not be smooth in space but analytic in time.
	Finally, Section \ref{4thsec} is about the time analyticity of some nonlinear heat equations with power nonlinearity of order $p$ \eqref{equa0}. We prove that a solution $u=u(x,t)$ of \eqref{equa0} is analytic in time $t\in (0,1]$ if it is bounded in $\mathrm{M}\times[0,1]$ and $p$ is a positive integer. In addition, we investigate the case when $p$ is a rational number with a stronger assumption $0<C_3 \leq |u(x,t)| \leq C_4$. As necessary lemmas, for any nonnegative integer $k$, we establish an explicit estimate of $|\p_t^k \Gamma(x,t;y)|$ where $\Gamma(x,t;y)$ is the heat kernel on a manifold, and a connection between $\p_t^k(t^k u^p)$ and $\p_t^k(t^k u)$.
	
	An interesting point is that the distribution of zeros of analytic solutions of the heat equation is connected to the Riemann Hypothesis as we can see from T.Tao$\&$B.Rodgers\cite{[Tao]} and V.G. Papanicolaou, E.Kallitsi$\&$G.Smyrlis\cite{[Vas]}.
	
	For the notation of this paper, we use $B(x, r)$ to denote the geodesic ball of radius $r$ centered at $x$ and
	$|B(x, r)|$ to denote the volume. $d(x, y)$ means the geodesic distance of
	$x$, $y$ $\in\mathrm{M}$ and $0$ denotes a reference point in $\mathrm{M}$. Besides, $Q_r(x,t)=B(x,r)\times(t-r^2,t)$ and $Q'_r (x,t)=B(x,r)\times(t-r^4,t)$.  
	Please note throughout this paper, constant $C$ may be different from case to case.

	\section{Biharmonic heat equation}\label{2ndsec}

	We now begin investigating the time analyticity of the biharmonic heat equation \eqref{2ndeq}. The main result in this section is Theorem \ref{2ndtheo}. First, we have several remarks about Theorem \ref{2ndtheo}.
	\begin{remark}
		Just note we use the condition that $u$ is of exponential growth of order $\frac{4}{3}$ in the computation of $\iint_{\Gamma_k^1}\left(u(x,t)\right)^2dxdt$ in \eqref{com}.
	\end{remark}
	\begin{remark}
		For any smooth solution $u=u(x,t)$ of the biharmonic heat equation \eqref{2ndeq} and any $(x_0,t_0)\in \mathrm{M}\times (0,1]$, actually we can get 
		$$|\partial_t^ku(x_0,t_0)| \leq \frac{A^*A_3^{k+1}k^k}{t_0^{k+q/4-d/8}}e^{2A_2d^{4/3}(x_0,0)},$$
		where $q=\left[\frac{d}{2}\right]+1$ and $[\cdot]$ means the floor function.
		Thus, we can see at $t=0$, this method fails to prove the time analyticity.
	\end{remark}
	\begin{remark}
		Just note the radius of convergence does not depends on $x$ because $A_3$ is independent of $x$. 
	\end{remark}
	\begin{remark}
		The exponential growth of order $\frac{4}{3}$ corresponds to the heat kernel estimate of the biharmonic heat equation \eqref{2ndeq} which can be found in G.Barbatis$\&$E.Davies\cite{[BD]}. Actually, we can expect that the solutions of high order Laplacian heat equation $u_t+(-\Delta)^m u=0$ are also analytic in time with exponential growth of order $\frac{2m}{2m-1}$ for any integer $m\geq 1$.
	\end{remark}
	\begin{remark}\label{rem}
		Now we briefly go over the main idea of the proof of Theorem \ref{2ndtheo}. 
		For any $(x_0,t_0)\in \mathrm{M}\times (0,1]$ and positive integer $k$, consider some regions for any $j=1,2,\cdots,k$,\\
		$\Gamma_{j}^{1}=\left\{(x,t) | d(x,x_0)<\frac{j t_0^{1/4}}{(2k)^{1/4}}, t \in[t_0-\frac{jt_0}{2k}, t_0]\right\},$\\
		$\Gamma_{j}^{2}=\left\{(x,t) | d(x,x_0)<\frac{(j+0.5)t_0^{1/4}}{(2k)^{1/4}},t \in[t_0-\frac{(j+0.5)t_0}{2k}, t_0]\right\}.$\\
		Immediately $\Gamma_{j}^{1}\subset \Gamma_{j}^{2} \subset \Gamma_{j+1}^{1}$.
		
		There are three main steps. We have a lemma for each step in the following.\\
		The first step is to prove that for some constant $C=C(d,K_0)$ and any $j=1,2,\cdots,k$,
		$$\iint_{\Gamma_{j}^{1}}|u_t(x,t)|^2dxdt\leq \frac{Ck}{t_0} \iint_{\Gamma_{j}^{2}} |\Delta u(x,t)|^2dxdt.$$
		The second step is to prove 
		$$\iint_{\Gamma_{j}^{2}} |\Delta u(x,t)|^2dxdt\leq \frac{Ck}{t_0} \iint_{\Gamma_{j+1}^{1}} |u(x,t)|^2dxdt.$$
		Then we can combine the above two inequalities
		and iterate to deduce
		$$\iint_{\Gamma_{1}^{1}}|\p_t^k u(x,t)|^2dxdt\leq \left(\frac{Ck}{t_0}\right)^{2k} \iint_{\Gamma_{k+1}^{1}} |u(x,t)|^2dxdt.$$\\
		The last step is to use the mean value inequality to get, for some constant $C=C(d,x_0,\mathrm{M})$,
		$$|\p_t^k u(x_0,t_0)|^2\leq C\left(\frac{k}{t_0}\right)^{1+q/2}\iint_{\Gamma_{1}^{1}}|\p_t^l u(x,t)|^2dxdt\leq C\left(\frac{k}{t_0}\right)^{1+q/2}\left(\frac{Ck}{t_0}\right)^{2k} \iint_{\Gamma_{k+1}^{1}} |u(x,t)|^2dxdt,$$
		which is exactly what we want.
	\end{remark}

	\subsection{Iterated energy estimates}
	Now we begin to estimate the $L_{loc}^2$ norm of $|\p_tu(x,t)|^2$.
	\begin{lemma}\label{1stlemmaine}
		For any smooth solution $u=u(x,t)$ of the biharmonic heat equation \eqref{2ndeq} and any $l=1,2,\cdots,k$, there exist some constant $C$ such that
		$$
		\iint_{\Gamma_{j}^{1}}|\p_tu(x,t)|^2dxdt\leq \frac{Ck}{t_0} \iint_{\Gamma_{j}^{2}}|\Delta u(x,t)|^2dxdt.
		$$
	\end{lemma}
	
	\pf By Theorem 6.33 of the paper J.Cheeger$\&$T. H.Colding\cite{[CC]}, there exists some smooth cut-off function $\psi^{(1)}(x,t)$ such that for some constant $C$,\\
	\be\label{condi1}
	\frac{|\nabla\psi^{(1)}(x,t)|^2}{\psi^{(1)}(x,t)}\leq \frac{C\sqrt{k}}{\sqrt{t_0}},\quad |\partial_t\psi^{(1)}(x,t)|+|\nabla \psi^{(1)}(x,t)|^4+|\Delta \psi^{(1)}(x,t)|^2\leq \frac{Ck}{t_0},
	\ee
	and 
	$$
	\begin{aligned}
	&0\leq\psi^{(1)}(x,t)\leq 1,\ \psi^{(1)}(x,t)=1 \ \text{in}\  \Gamma_{j}^{1},\ \psi^{(1)}(x,t) \ \text{is supported in}\ \Gamma_{j}^{2}.
	\end{aligned}
	$$
	As we are doing the biharmonic heat equation instead of the heat equation, we need to have the estimate for $|\Delta \psi^{(1)}(x,t)|^2$ which is why we need to cite the paper J.Cheeger$\&$T. H.Colding\cite{[CC]}.\\
	We use $\psi$ instead of $\psi^{(1)}(x,t)$ in this proof for simplicity of notation. By Green's formula, integration by parts and equation \eqref{2ndeq}, we find

	\be\label{tem}
	\begin{aligned}
		&\quad\iint_{\Gamma_{j}^{2}}|\p_tu(x,t)|^2\psi^2dxdt=-\iint_{\Gamma_{j}^{2}}\p_tu(x,t)\Delta^2u(x,t)\psi^2dxdt\\
		&=-\iint_{\Gamma_{j}^{2}}\Delta u(x,t)\Delta(\p_tu(x,t)\psi^2)dxdt\\
		&=-\iint_{\Gamma_{j}^{2}}\Delta u(x,t)\left(\Delta \p_tu(x,t)\psi^2+2\nabla \p_tu(x,t)\nabla \psi^2+\p_tu(x,t)\Delta\psi^2\right)dxdt\\
		&=-\frac{1}{2}\iint_{\Gamma_{j}^{2}}\p_t(\Delta u(x,t))^2\psi^2dxdt-2\iint_{\Gamma_{j}^{2}}\Delta u(x,t)\nabla \p_tu(x,t)\nabla\psi^2dxdt\\
		&\quad -\iint_{\Gamma_{j}^{2}}\Delta u(x,t)\p_tu(x,t)\Delta\psi^2dxdt\\
		&=\frac{1}{2}\iint_{\Gamma_{j}^{2}}(\Delta u(x,t))^2\p_t\psi^2dxdt-\frac{1}{2}\int_{B(x_0,\frac{(j+0.5)t_0^{1/4}}{(2k)^{1/4}})}(\Delta u(x,t))^2dx\bigg|_{t=t_0}
		\\
		&\quad -2\iint_{\Gamma_{j}^{2}}\Delta u(x,t)\nabla \p_tu(x,t)\nabla \psi^2dxdt\\
		&\quad-2\iint_{\Gamma_{j}^{2}}\Delta u(x,t)\p_tu(x,t)\Delta\psi \psi dxdt-2\iint_{\Gamma_{j}^{2}}\Delta u(x,t)\p_tu(x,t)|\nabla \psi|^2dxdt\\
		&\leq \frac{1}{2}\iint_{\Gamma_{j}^{2}}(\Delta u(x,t))^2\p_t\psi^2dxdt
		+2\iint_{\Gamma_{j}^{2}} \p_tu(x,t)\nabla \Delta u(x,t)\nabla \psi^2dxdt\\
		&\quad+2\iint_{\Gamma_{j}^{2}}\Delta u(x,t)\p_tu(x,t)\Delta\psi \psi dxdt+2\iint_{\Gamma_{j}^{2}}\Delta u(x,t)\p_tu(x,t)|\nabla \psi|^2dxdt.
	\end{aligned}
	\ee
	Next we can use the bounds for the cutoff function $\psi$ and the Cauchy-Schwarz inequality to get:
	\be\label{imp11}
	\begin{aligned}
		&\quad\iint_{\Gamma_{j}^{2}}|\p_tu(x,t)|^2\psi^2dxdt\\
		&\leq \frac{Ck}{t_0}\iint_{\Gamma_{j}^{2}} |\Delta u(x,t)|^2dxdt+\epsilon \iint_{\Gamma_{j}^{2}} |\p_tu(x,t)|^2\psi^2dxdt +\frac{4}{\epsilon}\iint_{\Gamma_{j}^{2}}|\nabla\Delta u(x,t)|^2|\nabla\psi|^2dxdt\\
		&+\epsilon\iint_{\Gamma_{j}^{2}}|\p_tu(x,t)|^2\psi^2dxdt+\frac{Ck}{\epsilon t_0}\iint_{\Gamma_{j}^{2}}|\Delta u(x,t)|^2dxdt\\
		&\quad+\epsilon\iint_{\Gamma_{j}^{2}}|\p_tu(x,t)|^2\psi^2 dxdt+\frac{Ck}{\epsilon t_0}\iint_{\Gamma_{j}^{2}}|\Delta u(x,t)^2|dxdt\\
		&= \frac{Ck}{t_0}(1+\frac{2}{\epsilon})\iint_{\Gamma_{j}^{2}} |\Delta u(x,t)|^2dxdt+3\epsilon\iint_{\Gamma_{j}^{2}}|\p_tu(x,t)|^2\psi^2dxdt\\
		&\quad+\frac{4}{\epsilon}\iint_{\Gamma_j^2}|\nabla\Delta u(x,t)|^2|\nabla\psi|^2dxdt.
	\end{aligned}
	\ee

	Now we need to get the estimate for the term $\frac{4}{\epsilon}\iint_{\Gamma_{j}^{2}}|\nabla\Delta u(x,t)|^2|\nabla\psi|^2dxdt$ as above. For some small positive constants $\epsilon_2$ and $\epsilon_3$,
	$$
	\begin{aligned}\label{imp12}
	&\quad \frac{4}{\epsilon}\iint_{\Gamma_{j}^{2}}|\nabla\Delta u(x,t)|^2|\nabla\psi|^2dxdt\leq \frac{4C\sqrt{k}}{\epsilon\sqrt{t_0}}\iint_{\Gamma_{j}^{2}}|\nabla\Delta u(x,t)|^2\psi dxdt\\
	&=-\frac{4C\sqrt{k}}{\epsilon\sqrt{t_0}}\iint_{\Gamma_{j}^{2}}\Delta^2u(x,t)\Delta u(x,t)\psi  dxdt\\
	&\quad-\frac{4C\sqrt{k}}{\epsilon\sqrt{t_0}}\iint_{\Gamma_{j}^{2}}\nabla\Delta u(x,t)\Delta u(x,t)\nabla\psi dxdt\\
	&\leq \frac{2C\sqrt{k}\epsilon_3}{\epsilon\sqrt{t_0}}\iint_{\Gamma_{j}^{2}}|\p_tu(x,t)|^2\psi^2dxdt+\frac{2C\sqrt{k}}{\epsilon\epsilon_3\sqrt{t_0}}\iint_{\Gamma_{j}^{2}}|\Delta u(x,t)|^2dxdt\\
	&\quad +\frac{2C\sqrt{k}\epsilon_2}{\epsilon\sqrt{t_0}}\iint_{\Gamma_{j}^{2}}|\nabla\Delta u(x,t)|^2|\nabla \psi|^2dxdt+\frac{2C\sqrt{k}}{\epsilon\epsilon_2\sqrt{t_0}}\iint_{\Gamma_{j}^{2}}|\Delta u(x,t)|^2dxdt.
	\end{aligned}
	$$
	Take $\epsilon=1/8$, $\epsilon_2=\frac{\sqrt{t_0}}{C\sqrt{k}}$ and $\epsilon_3=\frac{\sqrt{t_0}}{64C\sqrt{k}}$,
	we have       
	\be\label{imp13}
	\begin{aligned}
		&\frac{4}{\epsilon}\iint_{\Gamma_{j}^{2}}|\nabla\Delta u(x,t)|^2|\nabla\psi|^2dxdt\\
		&\leq\frac{1}{2}\iint_{\Gamma_{j}^{2}} |\p_tu(x,t)|^2\psi^2dxdt+\frac{2080C^2k}{t_0}\iint_{\Gamma_{j}^{2}}|\Delta u(x,t)|^2dxdt.
	\end{aligned}
	\ee
	By \eqref{imp11} and \eqref{imp13}, we can get
	$$\quad\iint_{\Gamma_{j}^{2}}|\p_tu(x,t)|^2\psi^2dxdt\leq \frac{Ck}{t_0}\iint_{\Gamma_{j}^{2}}|\Delta u(x,t)|^2dxdt,$$
	which finishes the proof of Lemma \eqref{1stlemmaine}.
	\qed
	
	Now we begin to estimate the $L_{loc}^2$ norm of $|\Delta u(x,t)|^2$. We can get a Caccioppoli type inequality (energy estimate) as follows.
	\begin{lemma}\label{2ndlemmaineq}
		For any  smooth solution $u=u(x,t)$ of the biharmonic heat equation \eqref{2ndeq} and any $l=1,2,\cdots,k$, there exist some constant $C$ such that
		\be\label{imp5}
		\sup_{t\in(t_0-\frac{(j+1)t_0}{(2k)},t_0)} \int_{B(x_0,\frac{(j+1)t_0^{1/4}}{(2k)^{1/4}})}u^2(x,t)\psi^2 dx+\iint_{\Gamma_{j}^{2}}|\Delta u(x,t)|^2dxdt\leq \frac{Ck}{t_0} \iint_{\Gamma_{j+1}^{1}}| u(x,t)|^2dxdt.
		\ee
	\end{lemma}
	
	\pf By Theorem 6.33 of the paper J.Cheeger$\&$T. H.Colding\cite{[CC]} again, there exists some smooth cut-off function $\psi^{(2)}(x,t)$ satisfying the condition \ref{condi1} and
	$$
	\begin{aligned}
	&0\leq\psi^{(2)}(x,t)\leq 1,\ \psi^{(2)}(x,t)=1 \ \text{in}\  \Gamma_{j}^{2},\  \psi^{(2)}(x,t) \ \text{is supported in}\ \Gamma_{j+1}^{1}.
	\end{aligned}
	$$
	
	We denote the cuf-off function $\psi^{(2)}(x,t)$ by $\psi$ again in this proof for the simplicity of notation. Similar  to \eqref{tem} and \eqref{imp11}, using Green's formula, Cauchy-Schwarz inequality, integration by parts and assumption for the cut-off function $\psi$, we yield
	\be\label{imp21}
	\begin{aligned}
		&\quad\iint_{\Gamma_{j+1}^{1}}(\Delta u(x,t))^2\psi^2dxdt\\
		&\leq \left(1+\frac{2}{\epsilon}\right)\frac{Ck}{t_0}\iint_{\Gamma_{j+1}^{1}}u^2(x,t)dxdt+3\epsilon \iint_{\Gamma_{j+1}^{1}}|\Delta u(x,t) |^2\psi^2dxdt\\
		&\quad+\frac{4}{\epsilon}\iint_{\Gamma_{j+1}^{1}}|\nabla u(x,t)|^2|\nabla \psi|^2dxdt,
	\end{aligned}
	\ee
	for any small positive constant $\epsilon$.
	
	Next we need to obtain the estimate for the term $\frac{4}{\epsilon}\iint_{\Gamma_{j+1}^{1}}|\nabla u(x,t)|^2|\nabla \psi|^2dxdt$.
	
	By integration by parts and Cauchy-Schwarz inequality, for some small positive constants $\epsilon_2$ and $\epsilon_3$,
	$$
	\begin{aligned}
	&\quad\frac{4}{\epsilon}\iint_{\Gamma_{j+1}^{1}}|\nabla u(x,t)|^2|\nabla \psi|^2dxdt\leq\frac{4C\sqrt{k}}{\epsilon\sqrt{t_0}}\iint_{\Gamma_{j+1}^{1}}|\nabla u(x,t)|^2 \psi dxdt\\
	&\leq \frac{2C\sqrt{k}\epsilon_2}{\epsilon\sqrt{t_0}}\iint_{\Gamma_{j+1}^{1}}|\Delta u(x,t)|^2\psi^2 dxdt+\frac{2C\sqrt{k}}{\epsilon\epsilon_2\sqrt{t_0}}\iint_{\Gamma_{j+1}^{1}}u^2(x,t) dxdt\\
	&\quad+\frac{2C\sqrt{k}\epsilon_3}{\epsilon\sqrt{t_0}}\iint_{\Gamma_{j+1}^{1}}|\nabla u(x,t)|^2|\nabla\psi|^2 dxdt+\frac{2C\sqrt{k}}{\epsilon\epsilon_3\sqrt{t_0}}\iint_{\Gamma_{j+1}^{1}}u^2(x,t)dxdt.
	\end{aligned}
	$$
	Take $\epsilon=\frac{1}{8}$, $\epsilon_2=\frac{\sqrt{t_0}}{128C\sqrt{k}}$ and $\epsilon_3=\frac{\sqrt{t_0}}{C\sqrt{k}}$, then
	\be\label{imp23}
	\begin{aligned}
		&\quad\frac{4}{\epsilon}\iint_{\Gamma_{j+1}^{1}}|\nabla u(x,t)|^2|\nabla \psi|^2dxdt\\
		&\leq  	\frac{1}{4}\iint_{\Gamma_{j+1}^{1}}|\Delta u(x,t)|^2\psi^2+\frac{4128C^2k}{t_0}\iint_{\Gamma_{j+1}^{1}}u^2(x,t)dxdt.
	\end{aligned}
	\ee
	Plugging \eqref{imp23} into \eqref{imp21}, we can get inequality 
	\be\label{imp2}
	\iint_{\Gamma_{j}^{2}}|\Delta u(x,t)|^2dxdt\leq \frac{Ck}{t_0} \iint_{\Gamma_{j+1}^{1}}| u(x,t)|^2dxdt.
	\ee
	
	Besides, we can also see 
	\be\label{notime}
	\begin{aligned}
		&\quad\p_t\left(1/2\int_{B(x_0,\frac{(j+1)t_0^{1/4}}{(2k)^{1/4}})}u^2(x,t)\psi^2 dx\right)\\
		&=\int_{B(x_0,\frac{(j+1)t_0^{1/4}}{(2k)^{1/4}})}-\Delta^2 u(x,t)u(x,t)\psi^2 dx+\int_{B(x_0,\frac{(j+1)t_0^{1/4}}{(2k)^{1/4}})}u^2(x,t)\psi\p_t\psi dx.\\
	\end{aligned}
	\ee
	For the term $\int_{B(x_0,\frac{(j+1)t_0^{1/4}}{(2k)^{1/4}})}-\Delta^2 u(x,t)u(x,t)\psi^2 dx$, we have, by integration by parts and assumption for $\psi$,
	$$
	\begin{aligned}
	&\quad\int_{B(x_0,\frac{(j+1)t_0^{1/4}}{(2k)^{1/4}})}-\Delta^2 u(x,t)u(x,t)\psi^2 dx\\
	&\leq \frac{1}{8}\int_{B(x_0,\frac{(j+1)t_0^{1/4}}{(2k)^{1/4}})}(\Delta(u(x,t)\psi))^2dx+\int_{B(x_0,\frac{(j+1)t_0^{1/4}}{(2k)^{1/4}})}|\nabla u(x,t)\nabla\psi|^2dx+\frac{Ck}{t_0}\int_{B(x_0,R)}|u(x,t)|^2 dx\\
	&\quad +\frac{1}{4}\int_{B(x_0,\frac{(j+1)t_0^{1/4}}{(2k)^{1/4}})} (u(x,t)\Delta\psi)^2dx+5\int_{B(x_0,\frac{(j+1)t_0^{1/4}}{(2k)^{1/4}})}(\Delta u(x,t)\psi)^2dx.
	\end{aligned}
	$$
	By integration about time in \eqref{notime}, using the assumption about $\psi$ and \eqref{imp23}, \eqref{imp2}, we can get the \eqref{imp5} immediately. 
	\qed
	
	\subsection{Mean value inequality for the biharmonic heat equation (\ref{2ndeq})}
	We also need the following lemma about the mean value inequality.
	\begin{lemma}
		Let $(x_0,t_0)$ be any point in $\mathrm{M}\times (0,1]$ and
		$u=u(x,t)$ be any solution to the biharmonic heat equation \eqref{2ndeq}. Then for some constant $C_1=C_1(d,x_0,\mathrm{M})$,
		\be\label{meann}
		\sup _{Q'_{r}\left(x_{0},t_0 \right)} |u(x,t)|^2 \leq \frac{C_1}{(R-r)^{2q+4}} \iint_{Q'_{R}\left(x_{0},t_0\right)} u^2(x,t) d x d t,
		\ee	
		where $q=\left[\frac{d}{2}\right]+1$ and $0<r<R<1$.
	\end{lemma}
	\begin{remark}
		Just note here the constant is dependent on $x_0$ and $\mathrm{M}$. This is because in the following proof, we need to use the Sobolev inequality, make sure the all the gradients of cut-off function $\psi$ below is bounded, and make sure $\nabla$ can commute with $\Delta$. In $R^d$, due to all of these peoperties satisfied, the constant $C$ should be independend of $x_0$ and $\mathrm{M}$.
	\end{remark}
	\pf

	Let $r<R_0<R_1<R_2<R$ where $R-R_2=R_2-R_1=R_1-R_0=R_0-r$ and define a smooth cut-off function $\phi=\phi(x,t)$ which is supported in $Q'_{R_0}(x_0,t_0)$ and $\phi=1$ in $Q'_{r}(x_0,t_0)$. Just note because the manifold is smooth in $B(x_0,1)$, for any nonnegative integer $k$, it holds for some constant $C=C(x_0,k,\mathrm{M})$,
	$$|\nabla^{k}R_m|\leq C(x_0,k,\mathrm{M}),$$
	where $R_m$ means the curvature tensor.
	
	Since $\phi$ is smooth in $B(x_0,1)$,  for any positive integer $i$, there exist some constant $C(x_0,i,\mathrm{M})$ depending on $x_0$, $i$ and $\mathrm{M}$ such that,
	$$|\nabla^i \phi^{2}|\leq \frac{C(x_0,i,\mathrm{M})}{|R_0-r|^i},\ |\Delta^i\phi|\leq \frac{C(x_0,2i,\mathrm{M})}{|R_0-r|^{2i}}$$
	$$|\nabla^i\p_t\phi|\leq \frac{C(x_0,4+i,\mathrm{M})}{|R_0-r|^{4+i}},\ \ |\Delta^i\p_t\phi|\leq \frac{C(x_0,2i+4,\mathrm{M})}{|R_0-r|^{2i+4}},$$
	where $\nabla$ is the covariant derivative and $\nabla^i$ means the i-th order covariant derivative.\\
	We can also define a smooth cut-off function $\psi=\psi(x,t)$ which is supported in $Q'_{R_1}(x_0,t_0)$ and $\psi=1$ in $Q'_{R_0}(x_0,t_0)$ satisfying similar condition as above.
	
	Following the method in H.Dong$\&$D.Kim\cite{[DK]}, we can use the Morrey type Sobolev inequality which can be find in Theorem 2.7 of E.Hebey\cite{[Hebey]}, which means there exist some constant $C=C(d,x_0,\mathrm{M})$ that
	$$\sup _{B(x_0,R_1)} |u(\cdot,t)\psi|\leq C\|u(\cdot,t)\psi\|_{W^{q,2}(B(x_0,R_1))}.$$
	
	Also, for some constant $C=C(d)$, by the fundamental theorem of calculus, we yield
	$$
	\begin{aligned}
	&\sup _{t\in(t_0-R_1^4,t_0)} |u(x,\cdot)\phi|^2\leq\sup _{t\in(t_0-R_1^4,t_0)}\int_{t_0-R_1^4}^{t} \p_t(u(x,\cdot)\phi)^2dt\\
	&\leq \int_{t_0-R_1^4}^{t_0}|\p_tu(x,t)|^2\phi^2dt+\frac{C}{(R-r)^4}\int_{t_0-R_1^4}^{t_0}|u(x,t)|^2dt
	\end{aligned}.$$
	Therefore for some $C=C(d,x_0,\mathrm{M})$,
	\be\label{add3rd}
	\sup _{Q'_{r}\left(x_{0},t_0 \right)} |u(x,t)|^2 \leq C\sum\limits_{i=0}^{q} \|\nabla^i(\p_tu(\cdot,\cdot)\psi+\frac{C}{(R-r)^4}u(\cdot,\cdot)\psi)\|_{W^{q,2}(Q_{R_1}(x_0,t_0))}.
	\ee
	Then we need to apply the well-known Bochner's formula and the related cummutation formula to commute $\nabla$ with $\Delta$ and its high order version, see Proposition 3.2.1 of Q.Zhang\cite{[Zhangbook]} e.g.. Using the above commutation formula,
	\be\label{add1st}
	\Delta\nabla_i f=\nabla_i\Delta f+R_{ij}\nabla_j f.
	\ee
	By this formula, for any smooth function $f$ and any cut-off function $\psi$ which is supported in $B(x_0,R_2)$, there exist some constant $C=C(d,K_0)$ such that
	\be\label{add2nd}
	\begin{aligned}
		&\int_{B(x_0,R_2)}(\Delta f)^{2}\psi^2 d x=\sum_{i, j=1}^{n} \int_{B(x_0,R_2)} \nabla_{i}\nabla_{i}f \nabla_{j}\nabla_{j}f\psi^2 d x\\
		&=-\sum_{i, j=1}^{n} \int_{B(x_0,R_2)} \nabla_{j}\nabla_{i}\nabla_{i}f\nabla_{j}f\psi^2 d x-\sum_{i, j=1}^{n} \int_{B(x_0,R_2)} \nabla_{i}\nabla_{i}f\nabla_{j}f\nabla_{j}\psi^2 d x\\
		&=\sum_{i, j=1}^{n} \int_{B(x_0,R_2)} \nabla_{i}\nabla_{j}f\nabla_{i}\nabla_{j}f\psi^2 d x+ \int_{B(x_0,R_2)}Ric\left(\nabla f,\nabla f\right)\psi^2dx\\
		&\quad+\sum_{i, j=1}^{n} \int_{B(x_0,R_2)} \nabla_{i}\nabla_{j}f\nabla_{j}f\nabla_{i}\psi^2 d x-\sum_{i, j=1}^{n} \int_{B(x_0,R_2)} \nabla_{i}\nabla_{i}f\nabla_{j}f\nabla_{j}\psi^2 d x\\
		&\geq 1/2\int_{B(x_0,R_2)}\left|\nabla^{2} f\right|^{2} \psi^2d x-C\int_{B(x_0,R_2)}\left|\nabla f\right|^2|\nabla\psi|^2dx.
	\end{aligned}
	\ee

	By using formula \eqref{add1st} and its high order version repeatedly and \eqref{add2nd} where we separate $(R_1,R_2)$ into $q$ equal parts, we can get for some $C=C(d,x_0,\mathrm{M})$,
	\be\label{add4th}
	\begin{aligned}
		&\sum\limits_{i=0}^{q}\iint_{Q'_{R_1}\left(x_{0},t_0\right)} |\nabla^i (\p_tu(x,t)\psi)|^2+ \frac{C}{(R-r)^4}|\nabla^i (u(x,t)\psi)|^2d x d t\\
		&\leq \sum\limits_{i=0}^{[\frac{q}{2}]}\frac{C}{(R-r)^{2q-4i}}\iint_{Q'_{R_2}\left(x_{0},t_0\right)} |\Delta^i \p_tu(x,t)|^2 d x d t\\
		&+\sum\limits_{i=0}^{[\frac{q}{2}]}\frac{C}{(R-r)^{2q-4i+4}}\iint_{Q'_{R_2}\left(x_{0},t_0\right)} |\Delta^i u(x,t)|^2 d x d t\\
		&+ \sum\limits_{i=0}^{[\frac{q-1}{2}]}\frac{C}{(R-r)^{2q-4i-1}}\iint_{Q'_{R_2}\left(x_{0},t_0\right)} |\nabla\Delta^i \p_tu(x,t)|^2 d x d t\\
		&+ \sum\limits_{i=0}^{[\frac{q-1}{2}]}\frac{C}{(R-r)^{2q-4i+3}}\iint_{Q'_{R_2}\left(x_{0},t_0\right)} |\nabla\Delta^i u(x,t)|^2 d x d t.
	\end{aligned}
	\ee
	By Lemma \ref{1stlemmaine} and Lemma \ref{2ndlemmaineq}, we have for some constant $C$
	$$\iint_{Q'_{R_2}\left(x_{0},t_0\right)} |u_t(x,t)|^2 d x d t\leq \frac{C}{|R-r|^4}\iint_{Q'_{R_3}\left(x_{0},t_0\right)} |u(x,t)|^2 d x d t.$$
	Plugging \eqref{imp2} and \eqref{imp23} into \eqref{add4th}, it holds for some constant $C=C(x_0,\mathrm{M})$, 
	$$\sum\limits_{i=0}^{q}\iint_{Q'_{R_1}\left(x_{0},t_0\right)} |\nabla^i (\p_t(u(x,t)\psi)+u(x,t)\psi)|^2 d x d t\leq \frac{C}{|R-r|^{4+2q}}\iint_{Q'_{R_3}\left(x_{0},t_0\right)} |u(x,t)|^2 d x d t.$$
	Plugging into \eqref{add3rd}, we can get \eqref{meann} immediately.

	\subsection{Proof of Theorem \ref{2ndtheo}}
	Now we are ready to prove Theorem \ref{2ndtheo}.
	Combining Lemma \ref{1stlemmaine} and Lemma \ref{2ndlemmaineq}, for any $l=1,2,\cdots,k$, we yield
	$$
	\iint_{\Gamma_j^1}|\p_tu(x,t)|^2dxdt\leq \frac{C^2k^2}{t_0^2} \iint_{\Gamma_{j+1}^1}| u(x,t)|^2dxdt.
	$$
	Since $\partial_t^l u $  is also a solution of \eqref{2ndeq}, 
	by iteration, we have
	$$
	\iint_{\Gamma_{1}^{1}}\left(\partial_t^k u(x,t)\right)^2 d x d t\leq \frac{C^2k^2}{t_0^2}\iint_{\Gamma_{2}^{1}}\left(\partial_t^{k-1} u(x,t)\right)^2 d x d t\leq...\leq \left(\frac{C^2k^2}{t_0^2}\right)^{k}\iint_{\Gamma_{k+1}^{1}} u(x,t)^2 d x d t.
	$$
	Using the mean value inequality \eqref{meann}, for some constant $A_3=A_3(d,K_0,A_2)$ and $A^*=A^*(A_1,d,x_0,\mathrm{M})$,
	\be\label{com}
	\begin{aligned}
		&|\partial_t^ku(x_0,t_0)|^2\leq C_1\left(\frac{2k}{t_0}\right)^{\frac{4+2q}{4}}  \iint_{Q'_{\left(\frac{t_0}{2k}\right)^{1/4}}\left(x_{0},t_0\right)} |\partial_t^ku(x,t)|^2 d x d t\\
		&\leq C_1\left(\frac{2k}{t_0}\right)^{\frac{4+2q}{4}}\left(\frac{C^2k^2}{t_0^2}\right)^{k}\iint_{\Gamma_k^1}\left(u(x,t)\right)^2dxdt\\
		&\leq C_1 \left(\frac{2k}{t_0}\right)^{\frac{4+2q}{4}}\left(\frac{C^2k^2}{t_0^2}\right)^{k} \times A_1^2 e^{4A_2d^{4/3}(x_0,0)}e^kt_0^{1+\frac{d}{4}}\leq \frac{{A^*}^2A_3^{2k+2}k^{2k}}{t_0^{2k+q/2-d/4}}e^{4A_2d^{4/3}(x_0,0)}.
	\end{aligned}
	\ee
	
	Thus,
	\begin{equation}\label{last1}
	|\partial_t^ku(x_0,t_0)| \leq \frac{A^*A_3^{k+1}k^{k}}{t_0^{k+q/4-d/8}}e^{2A_2d^{4/3}(x_0,0)}.
	\end{equation}
	Then we fix a number $R \geq 1$ and let $t \in[1-\delta, 1]$ for some small $\delta>0$.  For any positive integer $j$, Taylor's theorem implies that
	\begin{equation}\label{taylor}
	u(x,t)-\sum_{i=0}^{j-1} \partial_{t}^i u(x,1) \frac{(t-1)^{i}}{i !}=\frac{(t-1)^{j}}{j !} \partial_{t}^{j} u(x,s),
	\end{equation}
	where $s=s(x, t, j) \in[t, 1]$. By \eqref{last1}, for sufficiently small $\delta>0,$ the right-hand side of
	\eqref{taylor} converges to 0 uniformly for $x \in B(0, R)$ as $j \rightarrow \infty$. Hence
	$$
	u(x,t)=\sum_{j=0}^{\infty} \partial_{t}^{j} u(x,1) \frac{(t-1)^{j}}{j !}
	$$
	i.e., $u$ is analytic in time with radius $\delta$. Denote $a_{j}=a_{j}(x)=\partial_{t}^{j} u(x,1) .$ By \eqref{last1} again, we
	have
	$$
	\partial_{t} u(x,t)=\sum_{j=0}^{\infty} a_{j+1}(x) \frac{(t-1)^{j}}{j !} \text { and } \Delta^2 u(x,t)=\sum_{j=0}^{\infty} \Delta^2 a_{j}(x) \frac{(t-1)^{j}}{j !}
	$$
	where both series converge uniformly for $(x,t) \in B(0, R)\times [1-\delta,1]$. Since $u$ is a solution of the biharmonic heat equation \eqref{2ndeq}, it implies $-\Delta^2 a_{j}(x)=a_{j+1}(x)$
	with
	$$
	\left|a_{j}(x)\right| \leq A_1A_3^{k+1}k^ke^{2A_2d^{4/3}(x,0)}.
	$$
	This completes the proof
	of Theorem \ref{2ndtheo}.
	\qed\\

	We can then reach two corollaries similar to Corollary 2.2 and Corollary 2.6 in the paper H.Dong$\&$Q.Zhang\cite{[Zhang]}.
	
	\begin{corollary}\label{cor1}
		The Cauchy problem for the backward biharmonic heat equation
		\begin{equation}\label{back1}
		\quad\left\{\begin{array}{l}
		\partial_{t} u-\Delta^2 u=0 \\
		u(x,0)=a(x)
		\end{array}\right.
		\end{equation}
		has a smooth solution of exponential growth of order $\frac{4}{3}$ in $ \mathrm{M}\times (0, \delta) $ for some $\delta>0$ if and only if for any integer $k\geq 0$,
		\begin{equation}\label{ifonly3}
		|\left(\Delta^2\right)^{k} a(x)| \leq A_3^{k+1}k^ke^{A_2d^{\frac{4}{3}}(x_0,0)}, \quad j=0,1,2, \ldots
		\end{equation}
		where $A_{2}, A_3$ are some positive constants.
	\end{corollary}
	
	\pf Suppose \eqref{back1} has a smooth solution of exponential growth of order $\frac{4}{3}$, say $u=u(x,t)$. Then $u(x,-t)$ is a solution of the biharmonic heat equation \eqref{2ndeq} with polynomial growth of order $\frac{4}{3}$. By Theorem \ref{2ndtheo},
	\eqref{ifonly3} follows as $(-1)^j(\Delta^2)^{j} a(x)=a_{j}(x)$ in the theorem.
	
	On the other hand, suppose \eqref{ifonly3} holds. Then it is easy to check that
	$$u(x,t)=\sum_{j=0}^{\infty} (-1)^j(\Delta^2)^{j} a(x) \frac{t^{j}}{j !}$$
	is a smooth solution of the biharmonic heat equation for $t \in[-\delta, 0]$ with $\delta$ sufficiently small. Indeed, the bounds \eqref{last1} guarantee that the above series and the series
	$$
	\sum_{j=0}^{\infty} (-1)^{j+1}(\Delta^2)^{j+1} a(x) \frac{t^{j}}{j !} \text { and } \sum_{j=0}^{\infty} (-1)^j(\Delta^2)^{j} a(x) \frac{\p_t t^{j}}{j !}
	$$
	all converge absolutely and uniformly in $B(0, R)\times [-\delta, 0]$ for any fixed $R>0 .$ Hence $\partial_{t} u+\Delta^2 u=0 .$ Moreover $u$ has exponential growth of order $\frac{4}{3}$ since
	$$
	|u(x,t)| \leq \sum_{j=0}^{\infty}\left|(\Delta^2)^{j} a(x)\right| \frac{t^{j}}{j !} \leq  \sum_{j=0}^{\infty} A_3^{j+1}j^je^{A_2d^{\frac{4}{3}}(x_0,0)}\frac{t^{j}}{j !}  \leq A_{3} e^{A_2d^{\frac{4}{3}}(x_0,0)}
	$$
	for some $A_3$
	provided that $t \in[-\delta, 0]$ with $\delta$ sufficiently small. Thus, $u(x,-t)$ is a solution to the Cauchy problem of the backward biharmonic heat equation \eqref{back1} of exponential growth of order $\frac{4}{3}$.
	\qed
	
	\begin{remark} 
		It is known that generally the Cauchy problem for the backward biharmonic heat equation is not solvable. We can expect this corollary can be used in control theory, Ricci flow, stochastic analysis and some other areas.	
	\end{remark}
	
	We have another corollary about time analyticity at initial time $t=0$.
	\begin{corollary}\label{cor2}
		For the Cauchy problem for the biharmonic heat equation
		\begin{equation}\label{Cauchy1}
		\quad\left\{\begin{array}{l}
		\partial_{t} u+\Delta^2 u=0 \\
		u(x,0)=a(x).
		\end{array}\right.
		\end{equation}
		
		It has a smooth solution $u=u(x,t)$ of exponential growth of order $\frac{4}{3}$, which is analytic in time in
		$ \mathrm{M}\times [0, \delta)$ for some $\delta>0$ with a radius of convergence independent of $x$ if and only if
		\begin{equation}\label{ifonly4}
		|\left(\Delta^2\right)^k a(x)| \leq A_3^{k+1}k^ke^{A_2d^{\frac{4}{3}}(x_0,0)}, \quad k=0,1,2, \ldots, \quad j=0,1,2, \ldots
		\end{equation}
		where  $A_{2}, A_3$ are some positive constants.
	\end{corollary}

	\pf Assuming ( \ref{ifonly4} ), it is well-known that the problem \eqref{Cauchy1} has a solution 
	$$u=u(x,t)=\int_{\mathrm{M}}p(x,t;y)a(y)dy,$$ 
	for some $\delta>0$ and $t\in[0,\delta]$ where $p(x,t;y)$ is the heat kernel for the biharmonic heat equation on $\mathrm{M}$. 
	
	By Corollary \ref{cor1}, the following backward problem also has a solution
	$$
	\left\{\begin{array}{l}
	\partial_{t} v-\Delta^2 v=0 \\
	v(x,0)=a(x)
	\end{array}\right.
	$$
	in $\mathrm{M}\times [0, \delta)$ for some sufficiently small $\delta>0$. Define the function $U=U(x,t)$ by
	$$
	U(x,t)=\left\{\begin{array}{ll}
	u(x,t), & t \in[0, \delta) \\
	v(x,-t), & t \in(-\delta, 0]
	\end{array}\right.
	$$
	It is straight forward to check that $U(x,t)$ is a solution of the biharmonic heat equation in $\mathrm{M}\times (-\delta,\delta)$.\\
	By the theorem \ref{2ndtheo}, $U(x,t)$ and hence $u(x,t)$ is analytic in time at $t=0$.
	
	On the other hand, suppose $u(x,t)$ is a solution of the equation \eqref{Cauchy1}, which is analytic in
	time at $t = 0$ with a radius of convergence independent of $x$. Then, by definition, $u$ has
	a power series expansion in a time interval $(-\delta, \delta)$, for some $\delta > 0$. Hence, \eqref{ifonly4} holds
	following the proof of Corollary \ref{cor1}.
	\qed
	
	\begin{remark}
		Recall the well-known Kovalevskaya counter-example
		$$
		\quad\left\{\begin{array}{l}
		\partial_{t} u-\Delta u=0,\quad\forall (x,t)\in R\times [0,1]\\
		u(x,0)=\frac{1}{1+x^2},
		\end{array}\right.
		$$
		which says there are no analytic solutions in a neighborhood of the origin. We can extend it to the case of the biharmonic heat equation.
		\begin{lemma}
			Any smooth solution to the biharmonic heat equation\ref{2ndeq}
			$$
			\quad\left\{\begin{array}{l}
			\partial_{t} u+\Delta^2 u=0,\quad\forall (x,t)\in R\times [0,1] \\
			u(x,0)=\frac{1}{1+x^2},
			\end{array}\right.
			$$
			is not analytic near origin.
		\end{lemma}
		Actually, if we have a analytic solution $u$ near original, we can define $$u(x,t)=\sum\limits_{k,l\geq 0}a_{kl}\frac{t^k}{k!}\frac{x^l}{l!}.$$
		By induction, we can prove 
		$$a(m,2n)=(-1)^{m+n}(4m+2n)!, \  \text{for any nonnegative integers}\ m,n.$$
		Therefore $$\frac{|a(n,4n)|}{n!(4n)!}=\frac{(8n)!}{n!(4n)!}\to \infty,$$
		e solution is not analytic near origin.\\
		This corollary partially solves the problem about time analyticity of the biharmonic heat equation at $t=0$. 	
	\end{remark}
	\begin{remark}
		We can give a non-uniqueness example similar to the well-known non-uniqueness example for heat equation by A.N.Tychonov. To be precise, when $\mathrm{M}=R^1$, we can give a solution $u$ of \eqref{2ndeq} which does not satisfy $|u(x,t)|\leq A_1e^{A_2|x|^{4/3}}$ in $R^1\times (0,1]$ and is not analytic at $t=0$. It is 
		$$u(x,t)=\sum\limits_{k=0}^\infty (-1)^kD_t^k g(t) \frac{x^{4k}}{(4k)!},$$
		where 
		$$
		g(t)=\left\{\begin{array}{ll}
		e^{-t^{-\alpha}}, & \text{for any}\quad \alpha>1,t>0\\
		0, & t\leq 0.
		\end{array}\right.
		$$
		We can prove for some positive constant $C$, 
		$$|D_t^k g(t)|\leq \frac{C^k k!}{t^k}e^{-\frac{1}{2t^\alpha}},$$
		and therefore by $k!/(4k)!\leq 1/(3k)!$
		$$|u(x,t)|\leq \sum\limits_{k=0}^\infty\frac{C^k k!}{t^k}e^{-\frac{1}{2t^\alpha}}\frac{x^{4k}}{(4k)!}\leq \sum\limits_{k=0}^\infty\frac{C^k}{t^k}e^{-\frac{1}{2t^\alpha}}\frac{x^{4k}}{(3k)!}\leq e^{\left(\frac{C x^{4}}{t}\right)^{1/3}-\frac{1}{2t^\alpha}}.$$
		This example also shows the non-uniqueness for \eqref{2ndeq} because obviously we have another solution $u=0$.
	\end{remark}
	\section{Heat equation with potentials}\label{3rdsec}
	In this section, we mainly investigate the time analyticity of the heat equation with potentials \eqref{3rdeq}. The main idea of this section is similar to the idea as explained in Remark \ref{rem} of Section 2. First, let us define the weak solution.
	
	\begin{definition}
		We say $u=u(x,t)\in L^2_{loc}((t_1,t_2), W^{1,2}_{loc}(\mathrm{M})) $ is a weak subsolution (weak supersolution) to \eqref{3rdeq} if it satisfies,	
		$$
		\begin{aligned}
		&-\int_{t_1}^{t_2}\int_{\mathrm{M}}u(x,t)\p_t\phi(x,t)dxdt+\int_{t_1}^{t_2}\int_{\mathrm{M}}\nabla u(x,t)\nabla\phi(x,t)dxdt\\
		&+\int_{t_1}^{t_2}\int_{\mathrm{M}}V(x)u(x,t)dxdt\leq 0\ (\geq 0), 
		\end{aligned}
		$$
		for any nonnegative $\phi \in C_c^\infty \left( \mathrm{M}\times (t_1,t_2) \right)$.
	\end{definition}
	Especially, if $\phi \in C_c^\infty \left( \mathrm{M}\times (t_1,t_2) \right)$ and $\phi(\cdot,\frac{t_1+t_2}{2})=1$, then we can prove
	$$
	\begin{aligned}
	&-\int_{t_1}^{\frac{t_1+t_2}{2}}\int_{\mathrm{M}}u(x,t)\p_t\phi(x,t)dxdt+\int_{\mathrm{M}}u(x,\frac{t_1+t_2}{2})dx+\int_{t_1}^{\frac{t_1+t_2}{2}}\int_{\mathrm{M}}\nabla u(x,t)\nabla\phi(x,t)dxdt\\
	&+\int_{t_1}^{\frac{t_1+t_2}{2}}\int_{\mathrm{M}}V(x)u(x,t)dxdt\leq 0\ (\geq 0)
	\end{aligned}
	$$
	by testing with $\phi\eta_j$ and taking the limit $j\to \infty$, where $\eta_j=\eta_j(t)\in C_c^\infty(t_1,t_2)$ is a  sequence of nonnegative functions satisfying 
	$$\lim\limits_{j\to \infty}\eta_j(t)= \chi_{(t_1,\frac{t_1+t_2}{2})}\ a.e.$$
	We say $u$ is a weak solution if it is both a weak subsolution and a weak supersolution.\\
	
	Now for Theorem \ref{3rdtheo}, we have some remarks first.
	
	\begin{remark}
		To be more precise, for any $(x_0,t_0)\in \mathrm{M}\times (0,1/2]$, then in Theorem \ref{3rdtheo}, it holds
		$$
		|\partial_t^ku(x_0,t_0)| \leq  \frac{B_1B_2^{k+1}k^{k}}{t_0^{k}}e^{B_3d^2(x_0,0)}.
		$$
		for some constants $B_1$, $B_2$ and $B_3$.
		Besides, in Theorem \ref{4ththeo}, it holds
		$$
		|\partial_t^ku(x_0,t_0)| \leq  \frac{B_1B_2^{k+1}k^{k}}{t_0^{k}}e^{B_3d^2(x_0,0)}.
		$$
		for some constants $B_1$, $B_2$ and $B_3$.
	\end{remark}
	\begin{remark}\label{remark}
		By the method of Steklov average, or to be more precise, by Theorem 4.1 of 
		Q.Hou$\&$L.Saloff-Coste\cite{[HL]} which states if the heat kernel $\Gamma_V$ of \eqref{3rdeq} satisfies the $L^2$
		Gaussian type upper bound and, for any weak solution $u$ of \eqref{3rdeq}, $\p_t^l u$ is also a weak solution of \eqref{3rdeq} for any $l=1,2,\cdots$. on the one hand, if $V\geq 0$ and if $\Gamma$ is the heat kernel of heat equation on the same manifold $\mathrm{M}$, then by maximal principle, $0\leq \Gamma_V\leq \Gamma$, which means $\Gamma_V$ satisfies this Gaussian type upper bound condition considering \eqref{im3} and the mean value inequality. On the other hand, if $V(\cdot)\in L^{q}(B(0,R^*))$ for some $q>\frac{d}{2}$, it is well known that $\Gamma_V$ also satisfies this Gaussian type upper bound condition. 
		Besides, we can prove $\p_t^l u\in L^2_{loc}(\mathrm{M}\times (0,1))$ by combining \eqref{3rdineq1} and \eqref{3rdineq2} next. Therefore, $\p_t^l u$ is locally H\"older continuous, which means $u$ is smooth in time. 
	\end{remark}
	Now we begin to investigate Theorem \ref{3rdtheo}.

	\begin{remark}\label{remark11}
		In Theorem \ref{3rdtheo}, we have an extra condition $\inf\limits_{x\in\mathrm{M}}|B(x,1)|>0$ to use the Sobolev inequality.
		To be precise,  by Theorem 3.6 of E.Hebey\cite{[Hebey]},
		we can see for any $\lambda\in(0,1)$, $q\geq 1$ and $\frac{1}{q}\leq \frac{1-\lambda}{d}$, there exists some constant $C=C(d,\mathrm{M})$ such that 
		$$\|u\|_{C^\lambda(\mathrm{M})}\leq C\|u\|_{W_{1,q}(\mathrm{M})}.$$
		Also by Proposition 3.7 of E.Hebey\cite{[Hebey]}, if $d>q\geq 1$, then for
		$q^*=\frac{qd}{d-q}$, we have	$$\|u\|_{L^{q^*}(\mathrm{M})}\leq C\|u\|_{W_{1,q}(\mathrm{M})}.$$	  
	\end{remark}

	Before embarking on the proof of theorem \eqref{3rdtheo}, we need to have some lemmas first.  The first one is about the Poincar\'e inequality which is a result of \cite{[B]} and we can find it in Theorem 5.6.5 of L.Saloff-Coste\cite{[Saloff]}, e.g.. 
	\begin{lemma}
		Let $\mathrm{M}$ be a manifold satisfying same conditions as above Theorem \ref{2ndtheo}.
		
		Then for any $1\leq p<\infty$, there exists some constant $C=C(d,p,K_0)$ such that for any ball $B(x_0,r)\subset \mathrm{M}$ where $0<r<4$, 
		\be\label{poincare}
		\int_{B(x_0,r)}|f(x)-f_{B(x_0,r)}|^p dx\leq C r^p \int_{B(x_0,r)}|\nabla f(x)|^p dx,
		\ee	
		where $f_{B(x_0,r)}=\frac{\int_{B(x_0,r)}f(x)dx}{|B(x_0,r)|}$ is the mean value of $f$ in $B(x_0,r)$.
	\end{lemma}
	
	Using this result, we have the following lemma about the Sobolev inequality:
	
	\begin{lemma}\label{extra1}
		Let $\mathrm{M}$ be a manifold satisfying the conditions as above Theorem \ref{2ndtheo}.
		
		Then for any $1\leq p<\infty$, $f\in C_c^\infty(B(x_0,r))$ where $B(x_0,r)\subset \mathrm{M}$ with $r\leq 1$, there exist some constants $\nu_p>p$ and $C=C(d,p,K_0)$ such that
		\be\label{sobolevfinal1}
		\left(\int_{B(x_0,r)} |f|^{\frac{p\nu_p}{\nu_p-p}} dx \right)^{\frac{\nu_p-p}{\nu_p}} \leq \frac{C r^p}{|B(x_0,r)|^{\frac{p}{\nu_p}}} \int_{B(x_0,r)}|\nabla f|^p dx.
		\ee	
	\end{lemma}
	\pf By Bishop-Gromov volume comparison theorem,
	\be\label{volume}
	|B(x,r)|\leq |B(x,s)|\left(\frac{r}{s}\right)^d\exp{\left((d-1)\sqrt{K_0}r\right)}\leq C|B(x,s)|\left(\frac{r}{s}\right)^d,
	\ee	
	when $0<s<r<4$.
	
	Combine \eqref{poincare} and \eqref{volume}, we can get \eqref{sobolevfinal1} by Theorem 5.2.6 of L.Saloff-Coste\cite{[Saloff]} immediately.
	\qed
	\begin{remark}
		Here $\nu_2=d$ when $d>2$ and $\nu_2$ can be some number which is close to $2$ when $d=1$ or $d=2$. We use this Sobolev inequality for the mean value inequality in Lemma \ref{lem3new}. Unlike Remark \ref{remark11}, this is true for all dimensions but with extra condition $r\leq 1$.  
	\end{remark}
	
	Now for any $(x_0,t_0)\in \mathrm{M}\times(0,1/2]$, we introduce some regions similar to F.Lin$\&$Q.Zhang\\\cite{[LZ]} first. For any positive integer $k$ and any $j=1,2,\cdots,k$,\\
	$H_{j}^{1}=\left\{(x,t) | d(x,x_0)<\frac{j\sqrt{t_0}}{\sqrt{2k}}, t \in[t_0-\frac{jt_0}{2k}, t_0+\frac{jt_0}{2k}]\right\},$\\
	$H_{j}^{2}=\left\{(x,t) | d(x,x_0)<\frac{(j+0.5)\sqrt{t_0}}{\sqrt{2k}},t \in[t_0-\frac{(j+0.5)t_0}{2k}, t_0+\frac{(j+0.5)t_0}{2k}]\right\}.$\\
	So immediately $H_{j}^{1}\subset H_{j}^{2} \subset H_{j+1}^{1}$.
	
	Then we have the following lemma to estimate $\iint_{H_{j}^1}|\p_t u(x,t)|^2 dxdt$.
	\begin{lemma}
		For any $j=1,2,\cdots,k$, there exists some positive constant $C$ such that
		\be\label{3rdineq1}
		\begin{aligned}
			&\iint_{H_{j}^1}|\p_t u(x,t)|^2 dxdt\\
			&\leq \frac{Ck}{t_0} \iint_{H_{j}^2}|\nabla  u(x,t)|^2dxdt
			+\frac{Ck}{t_0}\iint_{H_{j+1}^1}|V(x)|| u(x,t)|^{2}dxdt.
		\end{aligned}
		\ee
	\end{lemma}
	\pf Let us define a smooth cut-off function $\phi^{(1)}(x,t)$ such that $\phi^{(1)}(x,t)=1$ in $H_{j}^{1}$ and is supported in $H_{j}^{2}$. We can also suppose there is some constant $C$ such that 
	$$|\nabla\phi^{(1)}(x,t)|^2+|\p_t\phi^{(1)}(x,t)|\leq \frac{Ck}{t_0}.$$
	We use $\phi=\phi^{(1)}(x,t)$ below for the simplicity of notation. By assumption of cut-off function $\phi$ and Cauchy-Schwarz inequality, integration by parts in time,
	$$
	\begin{aligned}
	&\quad\iint_{H_{j}^{2}}|\p_tu(x,t)|^2\phi^2 dxdt=\iint_{H_{j}^2}\p_tu(x,t)(\Delta u(x,t)-V(x)u(x,t))\phi^2 dxdt\\
	&=-\iint_{H_{j}^2}\nabla \p_tu(x,t)\nabla u(x,t)\phi^2 dxdt  -\iint_{H_{j}^2} \p_tu(x,t)\nabla u(x,t)\nabla\phi^2 dxdt    \\
	&\quad-\iint_{H_{j}^2} V(x) \p_tu(x,t)u(x,t)\phi^2 dxdt\\
	&\leq \frac{3}{4}\iint_{H_{j}^2}|\p_tu(x,t)|^2\phi^2 dxdt+\frac{Ck}{t_0}\iint_{H_{j}^2}|\nabla u(x,t)|^2|dxdt\\
	&+\frac{Ck}{t_0}\iint_{H_{j}^2}|V(x)||u(x,t)|^2\phi dxdt.
	\end{aligned}
	$$
	
	Then we can get a Caccioppoli type inequality (energy estimate) as below.
	
	\begin{lemma}
		For any $j=1,2,\cdots,k$, there exists some positive constant $C$ such that
		\be\label{imp6}
		\begin{aligned}
			&\sup_{t\in(t_0-\frac{(j+1)t_0}{2k},t_0+\frac{(j+1)t_0}{2k})} \int_{B(x_0,\frac{(j+1)\sqrt{t_0}}{\sqrt{k}})}u^2(x,t)\phi^2 dx\\
			&\quad+\iint_{H_{j}^2}|\nabla u(x,t)|^2 dxdt+\iint_{H_{j}^2}|V(x)||u(x,t)|^2 dxdt\leq \frac{Ck}{t_0}\iint_{H_{j+1}^1}|  u(x,t)|^2dxdt\\
			&\quad+C \left({C^*}^{\frac{2q}{2q-d}}+{D^*}^{\frac{2}{2q-d}}{C^{**}}^{\frac{2(q-1)}{2q-d}}\left(d(x_0,0)+\frac{(j+0.5)\sqrt{t_0}}{\sqrt{k}}\right)^{2}\right)\iint_{H_{j+1}^1}| u(x,t)|^{2}dxdt.
		\end{aligned}
		\ee
	\end{lemma}
	\pf Let us define another smooth cut-off function $\phi^{(2)}(x,t)$ such that $\phi^{(2)}(x,t)=1$ in $H_{j}^{2}$ and is supported in $H_{j+1}^{1}$. We can also suppose there is some constant $C$ such that 
	$$|\nabla\phi^{(2)}(x,t)|^2+|\p_t\phi^{(2)}(x,t)|\leq \frac{Ck}{t_0}.$$
	We use $\phi=\phi^{(2)}(x,t)$ for the simplicity of notation in this proof. By integration by parts, assumption about $\phi$ and \eqref{3rdeq},
	\be\label{c0}
	\begin{aligned}
		&\quad\iint_{H_{j+1}^1}|\nabla u(x,t)|^2\phi^2 dxdt\leq \frac{1}{4}\iint_{H_{j+1}^1}|\nabla u(x,t)|^2\phi^2dxdt\\
		&+\frac{Ck}{t_0}\iint_{H_{j+1}^1}|u(x,t)|^2 dxdt+\iint_{H_{j+1}^1}|V(x)||u(x,t)|^2\phi^2 dxdt.
	\end{aligned}
	\ee
	
	Now we need to estimate the last term above. By H\"older inequality, interpolation inequality and Sobolev inequality, we know:
	
	\be\label{c1}
	\begin{aligned}
		&\int_{B(x_0,\frac{(j+0.5)\sqrt{t_0}}{\sqrt{k}})}V(x)|u(x,t)|^2\phi^2 dx\\
		&\leq \left(\int_{B(x_0,\frac{(j+0.5)\sqrt{t_0}}{\sqrt{k}})}|V(x)|^qdx\right)^{1/q}\left(\int_{B(x_0,\frac{(j+0.5)\sqrt{t_0}}{\sqrt{k}})}\left(|u(x,t)|^2\phi^2\right)^{\frac{q}{q-1}}dx\right)^{\frac{q-1}{q}}\\
		&\leq  C\left(\int_{B(x_0,\frac{(j+0.5)\sqrt{t_0}}{\sqrt{k}})}|V(x)|^qdx\right)^{1/q}\left(2\epsilon^2 \left(\int_{B(x_0,\frac{(j+0.5)\sqrt{t_0}}{\sqrt{k}})}|\nabla (\phi u(x,t))|^{2^*}dx\right)^{\frac{2}{2^*}}\right)\\
		&\quad+C\left(\int_{B(x_0,\frac{(j+0.5)\sqrt{t_0}}{\sqrt{k}})}|V(x)|^qdx\right)^{1/q}\left(C(d,q) \epsilon^{\frac{-2d}{2q-d}} \int_{B(x_0,\frac{(j+0.5)\sqrt{t_0}}{\sqrt{k}})}|\phi u(x,t)|^{2}dx\right)
		\\
		&\leq C\left(\int_{B(x_0,\frac{(j+0.5)\sqrt{t_0}}{\sqrt{k}})}|V(x)|^qdx\right)^{1/q} \\
		&\quad\times \left(2\epsilon^2 \left(\int_{B(x_0,\frac{(j+0.5)\sqrt{t_0}}{\sqrt{k}})}|\nabla (\phi u(x,t))|^2dx+\int_{B(x_0,\frac{(j+0.5)\sqrt{t_0}}{\sqrt{k}})}| (\phi u(x,t))|^2dx\right)\right)\\
		&\quad+C\left(\int_{B(x_0,\frac{(j+0.5)\sqrt{t_0}}{\sqrt{k}})}|V(x)|^qdx\right)^{1/q} \left(\epsilon^{\frac{-2d}{2q-d}} \int_{B(x_0,\frac{(j+0.5)\sqrt{t_0}}{\sqrt{k}})}|\phi u(x,t)|^{2}dx\right).
	\end{aligned}
	\ee
	
	By taking $\epsilon=\frac{1}{2\sqrt{C}} \left(\int_{B(x_0,\frac{(j+0.5)\sqrt{t_0}}{\sqrt{k}})}|V(x)|^qdx\right)^{-\frac{1}{2q}}$ and integrating with respect to time,
	\be\label{tnew}
	\begin{aligned}
		&\iint_{H_{j}^2}V(x)|u(x,t)|^2\phi^2 dxdt\leq  \frac{1}{2}\iint_{H_{j}^2}|\nabla  u(x,t)|^2\phi^2dxdt+\frac{Ck}{t_0}\iint_{H_{j}^2}|  u(x,t)|^2dxdt\\
		&\quad+C\left(\int_{B(x_0,\frac{(j+0.5)\sqrt{t_0}}{\sqrt{k}})}|V(x)|^qdx\right)^{\frac{2}{2q-d}} \iint_{H_{j}^2}|\phi u(x,t)|^{2}dxdt
		\\
		&\leq  \frac{1}{2}\iint_{H_{j}^2}|\nabla  u(x,t)|^2\phi^2dxdt+\frac{Ck}{t_0}\iint_{H_{j}^2}|  u(x,t)|^2dxdt\\
		&+\left({C^*}^{\frac{2q}{2q-d}}+C{D^*}^{\frac{2}{2q-d}}{C^{**}}^{\frac{2(q-1)}{2q-d}}\left(d(x_0,0)+\frac{(j+0.5)\sqrt{t_0}}{\sqrt{k}}\right)^{2}\right)\iint_{H_{j}^2}|\phi u(x,t)|^{2}dxdt.   
	\end{aligned}
	\ee

	Plugging into \eqref{c0}, we yield ,
	
	\be\label{3rdineq2}
	\begin{aligned}
		&\iint_{H_{j}^2}|\nabla u(x,t)|^2 dxdt\leq \frac{Ck}{t_0}\iint_{H_{j+1}^1}|  u(x,t)|^2dxdt\\
		&+C \left({C^*}^{\frac{2q}{2q-d}}+{D^*}^{\frac{2}{2q-d}}{C^{**}}^{\frac{2(q-1)}{2q-d}}\left(d(x_0,0)+\frac{(j+0.5)\sqrt{t_0}}{\sqrt{k}}\right)^{2}\right)\iint_{H_{j}^2}|\phi u(x,t)|^{2}dxdt.
	\end{aligned}
	\ee
	
	Besides, we can also see 
	\be\label{notime1}
	\begin{aligned}
		&\quad\p_t\left(1/2\int_{B(x_0,\frac{(j+1)\sqrt{t_0}}{\sqrt{k}})}u^2(x,t)\phi^2 dx\right)\\
		&=\int_{B(x_0,\frac{(j+1)\sqrt{t_0}}{\sqrt{k}})} u(x,t)(\Delta u(x,t)-V(x)u(x,t))\phi^2 dx+\int_{B(x_0,\frac{(j+1)\sqrt{t_0}}{\sqrt{k}})}u^2(x,t)\phi\p_t\phi dx\\
		&\leq \frac{1}{2}\int_{B(x_0,\frac{(j+1)\sqrt{t_0}}{\sqrt{k}})}|\nabla u(x,t)|^2\phi^2 dx+\frac{Ck}{t_0}\int_{B(x_0,\frac{(j+1)\sqrt{t_0}}{\sqrt{k}})}|u(x,t)|^2 dx
	\end{aligned}
	\ee
	By integration by time and \eqref{3rdineq2}, we can get the \eqref{imp6} immediately. 
	\qed

	Then we need the mean value inequality as follows.
	\begin{lemma}\label{lem3new}
		Assume $\mathrm{M}$ is a manifold satisfying same conditions as Theorem \ref{3rdtheo}. Let  $u=u(x,t)$ be a nonnegative weak subsolution to \eqref{3rdineq2}.  Then for any $0<p<\infty$, $0<r<R<1$ and $(x_0,t_0)\in \mathrm{M}\times(0,1/2]$, 
		$$\begin{aligned}
		&\sup _{Q_{r}\left(x_{0},t_0 \right)}|u(x,t)|^p\leq C\left(\frac{R^2}{|B(x_0,R)|^{\frac{2}{\nu_2}}}\right)^{\frac{1}{\theta^*-1}}\\
		& \times
		\left(  \frac{1}{|R-r|^2}+{C^*}^{\frac{2q}{2q-d}}+{D^*}^{\frac{2}{2q-d}}{C^{**}}^{\frac{2(q-1)}{2q-d}}\left(d(x_0,0)+\frac{(j+0.5)\sqrt{t_0}}{\sqrt{k}}\right)^{2}    \right)^{\frac{\theta^*}{\theta^*-1}}\\
		&\quad\times \iint_{Q_{R}\left(x_{0},t_0\right)} |u(x,t)|^p d x d t ,
		\end{aligned}$$
		where $\theta^*=1+\frac{2}{\nu_2}$. Here $\nu_2$ is defined in Lemma \ref{extra1}.
	\end{lemma}
	\pf We can prove this one by Moser iteration. By H\"{o}lder inequality and Lemma \ref{extra1}, we have for any $w(x)\in C_c^\infty(B(x_0,R))$,
	\be\label{sobolevineq32}
	\begin{aligned}
		&\int_{B(x_0,R)}|w(x)|^{2(1+\frac{2}{\nu_2})}dx\leq \left(\int_{B(x_0,R)}|w(x)|^{\frac{2\nu_2}{\nu_2-2}}dx\right)^{\frac{\nu_2-2}{\nu_2}}\left(\int_{B(x_0,R)}|w(x)|^2dx\right)^{\frac{2}{\nu_2}}\\
		&\leq C\left(\frac{R^2}{|B(x_0,R)|^{\frac{2}{\nu_2}}}\right)\int_{B(x_0,R)}|\nabla w(x)|^2 dx\left(\int_{B(x_0,R)}|w(x)|^2dx\right)^{\frac{2}{\nu_2}}
	\end{aligned}
	\ee
	
	Let $\psi=\psi(x,t)$ be a standard smooth cut-off function such that $\psi=1$ in $Q_r(x_0,t_0)$ and is supported in $Q_R(x_0,t_0)$.
	We can assume $|\nabla\psi|^2+|\p_t\psi|\leq \frac{C}{|R-r|^2}$. Then by integration by parts and assumption about $\psi$ and \eqref{3rdeq},
	$$
	\begin{aligned}
	&\quad\p_t\left(1/2\int_{B(x_0,R)}|u(x,t)|^2\psi^2 dx\right)+\int_{B(x_0,R)}|\nabla(u(x,t)\psi)|^2dx\\
	&\leq \frac{C}{|R-r|^2}\int_{B(x_0,R)}|u(x,t)|^2dx+4\int_{B(x_0,R)}|\nabla u(x,t)|^2\psi^2 dx+\int_{B(x_0,R)} |V(x)||u(x,t)|^2\psi^2dx .
	\end{aligned}
	$$

	Combining \eqref{c1} and \eqref{3rdineq2}, by integrating with respect to time, we yield
	$$
	\begin{aligned}
	&\sup_{t\in(-R^2,0)}\int_{B(x_0,R)}|u(x,t)|^2\psi^2 dx+ \iint_{Q_R(x_0,t_0)}|\nabla(u(x,t)\psi)|^2dxdt\\
	&\leq C \left(\frac{1}{|R-r|^2}+{C^*}^{\frac{2q}{2q-d}}+{D^*}^{\frac{2}{2q-d}}{C^{**}}^{\frac{2(q-1)}{2q-d}}\left(d(x_0,0)+\frac{(j+0.5)\sqrt{t_0}}{\sqrt{k}}\right)^{2}\right)\\
	&\quad\times \iint_{Q_R(x_0,t_0)}|u(x,t)|^2dxdt.
	\end{aligned}
	$$
	
	Let $E(R)=\left(\frac{R^2}{|B(x_0,R)|^{\frac{2}{\nu}}}\right)$ and $F^*=\frac{1}{|R-r|^2}+{C^*}^{\frac{2q}{2q-d}}+{D^*}^{\frac{2}{2q-d}}{C^{**}}^{\frac{2(q-1)}{2q-d}}\left(d(x_0,0)+\frac{(j+0.5)\sqrt{t_0}}{\sqrt{k}}\right)^{2}$ for simplicity of notataion.
	From inequality \eqref{sobolevineq32}, we can see 
	\be\label{in}
	\begin{aligned}
		&\iint_{Q_R(x_0,t_0)}|u(x,t)\psi|^{2\theta^*}dxdt\leq CE(R)\left(F^*\iint_{Q_R(x_0,t_0)}|u(x,t)|^2dxdt\right)^{\theta^*}.
	\end{aligned}
	\ee
	
	Now we have two cases.\\
	\textbf{Case (1):} $\boldsymbol{p \geq 2}$. In this case, we can see $u^{p/2}$ is also a nonnegative subsolution. Therefore,  \eqref{in} yields that:
	\be\label{in2}
	\iint_{Q'_r(x_0,t_0)}\left(u(x,t)\right)^{p\theta}dxdt
	\leq CE(R)\left(F^*\iint_{Q'_R(x_0,t_0)}(u(x,t) )^{p}dxdt\right)^{\theta}.
	\ee
	
	Set for some positive constant $\delta=\frac{r}{R}<1$, $\omega_{i}=\frac{(1-\delta)R}{2^{i}}$ so that $\sum_{1}^{\infty} \omega_{i}=(1-\delta)R$. Set also $\sigma_{0}=R, \sigma_{i+1}=$
	$\sigma_{i}-\omega_{i}=R-\sum_{1}^{i} \omega_{j} .$ Applying \eqref{in2} with $p=p_{i}=\theta^{i}, r=\sigma_{i}, R=\sigma_{i+1}$
	we obtain
	
	$$
	\iint_{Q'_{\sigma_{i+1}}(x_0,t_0)}\left(u(x,t)\right)^{\theta^{i+1}}dxdt
	\leq CE(R)16^{i\theta}{F^*}^\theta\left(\iint_{Q'_{\sigma_{i}}(x_0,t_0)}(u(x,t) )^{\theta^i}dxdt\right)^{\theta}.
	$$
	Hence, by iteration,
	$$
	\begin{aligned}
	&\quad\left(\iint_{Q'_{\sigma_{i+1}}(x_0,t_0)}\left(u(x,t)\right)^{\theta^{i+1}}dxdt\right)^{\theta^{-1-i}}\\
	&\leq (CE(R))^{\Sigma \theta^{-1-j}}16^{\Sigma j\theta^{-1-j}} {F^*}^{\Sigma \theta^{-j}} \iint_{Q'_{R}(x_0,t_0)}(u(x,t) )^{2}dxdt,
	\end{aligned}
	$$
	where all the summations are taken from 0 to $i$ and we can easily see $\Sigma j\theta^{-1-j}$ converges.  Letting $i$ tend to infinity, we obtain
	\be\label{mean1}
	\sup _{Q'_{r}(x_0,t_0)}u^{2} \leq C\left(E(R)\right)^{\frac{1}{\theta-1}}{F^*}
	^{\frac{\theta}{\theta-1}}\|u\|_{2, Q'_{R}(x_0,t_0)}^{2}.
	\ee
	Then when $p>2$, we can see $u^{p/2}$ is also a nonnegative subsolution, so 
	$$
	\sup_{Q'_{r}(x_0,t_0)}u^{p} \leq C\left(E(R)\right)^{\frac{1}{\theta-1}}{F^*}
	^{\frac{\theta}{\theta-1}}\|u\|_{p, Q'_{R}(x_0,t_0)}^{p},
	$$
	which proves \eqref{meann} for the case $p\geq 2$.

	\textbf{Case (2):} $\boldsymbol{0<p<2}$.\\
	For this case, we can use the method of M.Giaquinta\cite{[Giaquinta]} or more precisely, Theorem 2.2.3 in the book L.Saloff-Coste\cite{[Saloff]}.
	
	Fix $\sigma \in(0,1)$ and set $\rho=\sigma+(1-\sigma) / 4.$ Then \eqref{mean1} applies
	$$
	\sup _{Q'_{\sigma R}(x_0,t_0)}u \leq C\left(E(R)\right)^{\frac{1}{2(\theta-1)}}\left(1+\frac{1}{(\rho R-\sigma R)^4}\right)^{\frac{\theta}{2\theta-2}}\|u\|_{2, Q'_{\rho R}(x_0,t_0)}.
	$$
	Now, as $\|u\|_{2, Q} \leq\|u\|_{\infty, Q}^{1-p / 2}\|u\|_{p, Q}^{p / 2}$ for any parabolic cylinder Q, we
	get
	\be\label{ite}
	\|u\|_{\infty, Q'_{\sigma R}(x_0,t_0)} \leq J\left(1+\frac{1}{(\rho R-\sigma R)^4}\right)^{\frac{\theta}{2\theta-2}}\|u\|_{\infty, Q'_{\rho R}(x_0,t_0)}^{1-p / 2},
	\ee
	where $J=C\|u\|_{p, Q'_{\rho R}(x_0,t_0)}^{p / 2}\left(E(R)\right)^{\frac{1}{2(\theta-1)}}$.
	
	Fix $\delta =\frac{r}{R}$, $\sigma_{0}=\delta R=r$ and $\sigma_{i+1}=\sigma_{i}+\left(R-\sigma_{i}\right) / 4 .$ 
	Then
	$R-\sigma_{i}=$
	$(3 / 4)^{i}(1-\delta)R .$\\
	
	Applying the above inequality \eqref{ite} for each $i$ yields
	$$
	\|u\|_{\infty,  Q'_{\sigma_i}(x_0,t_0)} \leq(4 / 3)^{\theta i / (2\theta-2)} J{F^*}^{\frac{\theta}{2\theta-2}}\|u\|_{\infty,  Q'_{\sigma_{i+1}}(x_0,t_0)}^{1-p / 2}.
	$$
	Hence by iteration, for $i=1,2, \dots$
	$$
	\|u\|_{\infty, Q'_{r}(x_0,t_0)} \leq(4 / 3)^{(\theta / (\theta-2)) \sum_{0}^{i-1} j(1-p / 2)^{j}}\left[J{F^*}^{\frac{\theta}{2\theta-2}}\right]^{ \sum_{0}^{i-1}(1-p / 2)^{j}}\|u\|_{\infty,  Q'_{\sigma_i}(x_0,t_0)}^{(1-p / 2)^{i}}.
	$$
	Letting $i$ tend to infinity, we yield,
	$$
	\|u\|_{\infty, Q'_{r}(x_0,t_0)} \leq C\left(E(R)\right)^{\frac{1}{p(\theta-1)}}{F^*}^{\frac{\theta}{(\theta-1)p}}\|u\|_{p, Q'_{R}(x_0,t_0)},
	$$
	which proves inequality \eqref{meann} for the case $0<p<2$.
	
	\qed

	\subsection{Proof of Theorem \ref{3rdtheo}} Now we are in a position to prove Theorem \ref{3rdtheo}.
	
	Because $\p_t^l u(x,t)$ is also a weak solution of \eqref{3rdeq} for any $l=1,2,\cdots,k$, we can put inequality \eqref{3rdineq1} and \eqref{3rdineq2} together to obtain,
	$$
	\begin{aligned}
	&\iint_{H_{j}^1}|\p_t^{k-j+1} u(x,t)|^2 dxdt\leq \frac{C^2k^2}{t_0^{2k}} \iint_{H_{j+1}^1}| \p_t^{k-j} u(x,t)|^2dxdt\\
	&+\frac{Ck}{t_0}
	\left({C^*}^{\frac{2q}{2q-d}}+C{D^*}^{\frac{2}{2q-d}}{C^{**}}^{\frac{2(q-1)}{2q-d}}\left(d(x_0,0)+\frac{(j+0.5)\sqrt{t_0}}{\sqrt{k}}\right)^{2}\right)\iint_{H_{j+1}^1}| \p_t^{k-j} u(x,t)|^{2}dxdt\\
	&\leq \left(\frac{C^2k^2}{t_0^2}\right)\iint_{H_{j+1}^1}| \p_t^{k-j} u(x,t)|^{2}dxdt\\
	&+\left({C^*}^{\frac{4q}{2q-d}}+{D^*}^{\frac{4}{2q-d}}{C^{**}}^{\frac{4(q-1)}{2q-d}}\left(d(x_0,0)+\frac{(j+0.5)\sqrt{t_0}}{\sqrt{k}}\right)^{4}\right)\iint_{H_{j+1}^1}| \p_t^{k-j} u(x,t)|^{2}dxdt.
	\end{aligned}
	$$
	
	By iteration,
	\be\label{final}
	\begin{aligned}
		&\iint_{H_{1}^1}|\p_t^{k} u(x,t)|^2 dxdt\\
		&\leq \prod_{j=1}^{k}\left(\frac{C^2k^2}{t_0^2}+{C^*}^{\frac{4q}{2q-d}}+{D^*}^{\frac{4}{2q-d}}{C^{**}}^{\frac{4(q-1)}{2q-d}}\left(d(x_0,0)+\frac{(j+0.5)\sqrt{t_0}}{\sqrt{k}}\right)^{4}\right)\\
		&\quad\times\iint_{H_{k+1}^1}| u(x,t)|^{2}dxdt.
	\end{aligned}
	\ee
	By Lemma 5.2.7 of L.Saloff-Coste\cite{[Saloff]} or the book M.Richard$\&$S.T. Yau\cite{[SY]}, we see for some constant $D>0$ and any $0<r<1$, 
	\be\label{volume1}
	|B(x,r)|\leq e^{D \frac{d(x,y)}{r}}|B(y,r)|.
	\ee
	As $|\p_t^ku|^2$ is a weak subsolution to \eqref{3rdeq}, by mean value inequality in Lemma \ref{lem3new}, it holds
	\be\label{temp}
	\begin{aligned}
		&\quad|\p_t^ku(x_0,t_0)|^2\\
		& \leq Ce^{D d(x_0,0)}\left(\frac{k}{t_0}\right)^{\frac{d-\nu_2}{2}}\left(\frac{k}{t_0}+{C^*}^{\frac{2q}{2q-d}}+{D^*}^{\frac{2}{2q-d}}{C^{**}}^{\frac{2(q-1)}{2q-d}}\left(d(x_0,0)+\frac{(j+0.5)\sqrt{t_0}}{\sqrt{k}}\right)^{2}\right)^{\frac{\nu_2+2}{2}}\\
		&\quad\times\iint_{H_{1}^1}|\p_t^{k} u(x,t)|^2 dxdt.
	\end{aligned}
	\ee
	
	Combining these two inequalities \eqref{final} and \eqref{temp}, and applying the assumption that $u$ is of exponential growth of order $2$, we yield,
	$$
	|\partial_t^ku(x_0,t_0)|^2 \leq  \frac{A_1A_3^{2k+2}k^{2k}}{t_0^{2k}}e^{2A_4d^2(x_0,0)}.
	$$
	Just note we put some terms involving $d(x_0,0)$ into $e^{2A_4d^2(x_0,0)}$.
	
	The proof for the conclusions about $a_j=\p_t^ju(x,1/2)$ is the same as Theorem \ref{2ndtheo}. In this way, we have completed the proof of Theorem \ref{3rdtheo}.
	\qed
	
	\begin{remark}
		To see the set of functions satisfying the condition \ref{condnew} is nontrivial when $V(x)=x^2$ in $R^d$, we give some examples here.
		For Hermite polynomials $$H_n(x)=(-1)^n e^{x^2}\frac{d^n}{dx^n}e^{-x^2},$$
		and $$\psi_n(x)=(2^n n! \sqrt{\pi})^{-1/2}e^{-x^2/2}H_n(x),$$ it is well-known that
		$(D^2-x^2)\psi_n(x)=-(2n+1)\psi_n(x)$ and thus $$(D^2-x^2)^{k}\psi_n(x)=(-1)^k(2n+1)^k\psi_n(x).$$
		Therefore, $\psi_n(x)$ satisfies the condition \ref{condnew} as $|(D^2-x^2)^{k}\psi_n(x)|\leq C^k k!$.
	\end{remark}
	
	\begin{remark}
		Theorem \eqref{3rdtheo} is about the time analyticity when $t\in (0,1/2]$. Because \eqref{3rdeq} is a linear equation, it is a natural assumption that $u$ is of exponential growth of order 2 in $t\in[0,2]$, a longer time interval, then the solution should be time analytic in $[0,1]$.
	\end{remark}

	Especially, when $\mathrm{M}=R^d$, there is no necessity to assume $V\in L^1(R^d\ B(0,R^*))$ and $d\geq 3$, instead we have the following corollary:
	\begin{corollary}
		Let $\mathrm{M}=R^d$.
		Assume $V=V(x)$ satisfies the following conditions:\\
		(1) There exists some $R^*>0$ such that $V(\cdot)\in L^{q}(B(0,R^*))$ for some $q>\frac{d}{2}$. \\
		(2) For some constant $C^{**}>0$, if $d(x,0)>R^*$, then $|V(x)|\leq C^{**} d(x,0)^{\alpha}$ where $\alpha=2-\frac{2d}{q}$.\\
		Let $$\|V\|_{L^{q}(B(0,R^*))}=C^{*}$$ where $C^*$ is a positive constant and let $u=u(x,t)$ be a weak solution of equation \eqref{3rdeq} for any dimension $d\geq 1$ on $\mathrm{M}\times [0,1]$ of exponential growth of order $2,$ namely
		$$|u(x,t)| \leq A_{1} e^{A_{2} d^{2}(x, 0)}, \quad \forall(x,t) \in\mathrm{M}\times [0,1] , $$
		where $A_{1}$ and $A_{2}$ are some positive constants. Then $u$ is analytic in $t \in (0,1/2]$ with radius of convergence depending only on $t$, $d$, $q$, $K_0$, $A_{2}$, $\alpha$ and $C^*$.
		
		Moreover, if $t\in(1/2-\delta,1/2]$ for some small $\delta>0$, we have
		$$
		u(x,t)=\sum_{j=0}^{\infty} a_{j}(x) \frac{(t-1/2)^{j}}{j !}
		$$
		with $(\Delta-V) a_{j}(x)=a_{j+1}(x),$ and
		$$
		\left|a_{j}(x)\right| =\left|(\Delta-V)^j a_0(x)\right|\leq A_1A_{3}^{j+1} j^{j} e^{A_4 d^{2}(x, 0)}, \quad j=0,1,2, \ldots
		$$
		where constants $A_3=A_3(d,q,K_0,A_2,\alpha,C^*)$ and $A_4=A_4(A_2,\alpha,C^{**})$.\\
	\end{corollary}
	
	\pf\\ 
	The proof is almost the same as Theorem \ref{3rdtheo}. There are just two differences. The first one is to make a little change in
	\eqref{tnew}, instead, we yield,
	$$
	\begin{aligned}
	&C\left(\int_{B(x_0,\frac{(j+0.5)\sqrt{t_0}}{\sqrt{k}})}|V(x)|^qdx\right)^{\frac{2}{2q-d}} \iint_{H_{j}^2}|\phi u(x,t)|^{2}dxdt\\
	&\leq
	C\left({C^{**}}^{\frac{2q}{2q-d}}\left(\frac{(j+0.5)\sqrt{t_0}}{\sqrt{k}}+d(x_0,0)\right)^{2}+ {C^*}^{\frac{2q}{2q-d}}\right)\iint_{H_{j}^2}|\phi u(x,t)|^{2}dxdt.
	\end{aligned}
	$$
	The second difference is in \eqref{c1}. Instead of Sobolev inequality,  we use the Gagliardo-Nirenberg interpolation inequality and Young's inequality directly, which is
	$$
	\begin{aligned}
	&\|u(\cdot,t)\phi\|_{L^{\frac{2q}{q-1}}(R^d)}\leq C\|\nabla(u(\cdot,t)\phi)\|_{L^2(R^d)}^{\frac{d}{2q}}\|u(\cdot,t)\phi\|_{L^2(R^d)}^{1-\frac{d}{2q}}\\
	&\leq \epsilon\|\nabla(u(\cdot,t)\phi)\|_{L^2(R^d)}+C\epsilon^{\frac{-d}{2q-d}}\|u(\cdot,t)\phi\|_{L^2(R^d)}.
	\end{aligned}
	$$
	
	The rest of the proof is exact same.
	\qed\\

	As a special case, when $V(x)\geq 0$, we need to prove Theorem \ref{4ththeo} now.
	\begin{remark}
		In Theorem \ref{3rdtheo}, an interesting property is that the solution $u=u(x,t)$ can be not smooth in $x$ at all. Actually, if $\mathrm{M}=R^d$ and $V(x)=\frac{A}{|x|^2}$ where $A\geq 0$, we have one solution $u(x,t)=|x|^{\alpha(A)}$ where $\alpha(A):=\frac{-(d-2)+\sqrt{(d-2)^2+4A}}{2}$. We can see this solution is not smooth if $\alpha(A)$ is not an integer. 
	\end{remark}
	
	Similarly, we have a lemma about the mean value inequality using the same proof as in Lemma \ref{lem3new}:
	
	\begin{lemma}\label{lem4new}
		Assume $\mathrm{M}$ is a manifold satisfying same conditions as Theorem \ref{4ththeo}. Then
		for any nonnegative weak subsolution $u=u(x,t)$ to \eqref{3rdeq} where $V\geq 0$, for any $0<p<\infty$, $0<r<R<1$ and $(x_0,t_0)\in \mathrm{M}\times[-1,0]$, 
		there exist some canstant $C$ such that:
		$$\begin{aligned}
		&\sup _{Q_{r}\left(x_{0},t_0 \right)}|u(x,t)|^p \leq C\left(\frac{R^2}{|B(x_0,R)|^{\frac{2}{\nu}}}\right)^{\frac{1}{\theta^*-1}}
		\left(  \frac{1}{|R-r|^2}      \right)^{\frac{\theta^*}{\theta^*-1}}\iint_{Q_{R}\left(x_{0},t_0\right)} |u(x,t)|^p d x d t ,
		\end{aligned}$$
		where $\theta^*=1+\frac{2}{\nu_2}$ and $\nu_2$ is defined in Lemma \ref{extra1}. 
	\end{lemma}
	\begin{remark}
		As a very special example, we get the heat equation with inverse-square potential when $V(x)=\frac{A}{d(x,0)^2}$,
		$$
		\p_tu(x,t)-\Delta u(x,t)+\frac{Au(x,t)}{d(x,0)^2}=0, \quad\forall (x,t)\in \mathrm{M}\times [0,1].
		$$
		It is well-konwn that this potential is a borderline one where the regularity theory differs from the standard one. 
		For the regularity and mean value inequality of this equation in $R^d$, we can refer to Z.Li$\&$Q.Zhang\cite{[ZZ]}, B.Wong$\&$Q.Zhang\cite{[WZ]} and Z.Li$\&$X.Pan\cite{[LP]}. Actually, the inverse-square potential term $\frac{A}{|x|^2}$ helps with it.
	\end{remark}

	\subsection{Proof of Theorem \ref{4ththeo}}
	
	Now for any $(x_0,t_0)\in \mathrm{M}\times(0,1]$, we introduce some regions first. For any positive integer $k$ and any $j=1,2,\cdots,k$,\\
	$H_{j}^{1}=\left\{(x,t) | d(x,x_0)<\frac{j\sqrt{t_0}}{\sqrt{2k}}, t \in[t_0-\frac{jt_0}{2k}, t_0]\right\},$\\
	$H_{j}^{2}=\left\{(x,t) | d(x,x_0)<\frac{(j+0.5)\sqrt{t_0}}{\sqrt{2k}},t \in[t_0-\frac{(j+0.5)t_0}{2k}, t_0]\right\}.$\\
	So immediately $H_{j}^{1}\subset H_{j}^{2} \subset H_{j+1}^{1}$.
	
	Denote by $\psi_j^{(1)}(x,t)$ a standard smooth cut-off function supported in $H_{j}^{2}$
	such that\\  $\psi_j^{(1)}(x,t)=1$ in $H_{j}^{1}$ and $|\partial_t\psi_j^{(1)}(x,t)|+|\nabla \psi_j^{(1)}(x,t)|^2\leq \frac{Ck}{t_0}$ for some constant $C$.
	
	We denote $\psi=\psi_j^{(1)}(x,t)$ for simplicity of notation below. Then by equation \eqref{3rdeq} and integration by parts,
	
	$$
	\begin{aligned}
	& \iint_{H_{j}^{2}} (\p_t u(x,t))^2\psi^2 dxdt\leq \frac{1}{2}\iint_{H_{j}^{2}}|\nabla u(x,t)|\p_t \psi^2 dxdt\\
	& +\epsilon_1\iint_{H_{j}^{2}}(\p_t u(x,t))^2\psi^2 dxdt +\frac{4}{\epsilon_1}\iint_{H_{j}^{2}}|\nabla u(x,t)|^2 |\nabla \psi|^2 dxdt+\frac{1}{2}\iint_{H_{j}^{2}}V(x)u^2(x,t)\p_t\psi^2dxdt.
	\end{aligned}
	$$

	Using the assumption of $\psi$ and taking $\epsilon_1=\frac{1}{2}$, we yield 
	\be\label{imp41}
	\iint_{H_{j}^{1}}|\p_tu(x,t)|^2dxdt\leq \frac{Ck}{t_0} \left(\iint_{H_{j}^{2}}|\nabla u(x,t)|^2dxdt+\iint_{H_{j}^{2}}V(x) u^2(x,t)dxdt\right).
	\ee
	
	Define another smooth cut-off function $\psi_j^{(2)}(x,t)$ 
	supported in $H_{j+1}^{1}$ such that $\psi_j^{(2)}(x,t)=1$ in $H_{j}^{2}$. We assume for some constant $C$,  $|\partial_t\psi_j^{(2)}(x,t)|+|\nabla \psi_j^{(2)}(x,t)|^2\leq \frac{Ck}{t_0}$. We denote $\psi=\psi_j^{(2)}(x,t)$ for simplicity of notation below. Then by equation \eqref{3rdeq},
	$$
	\begin{aligned}
	&\quad \iint_{H_{j+1}^{1}}|\nabla u(x,t)|^2\psi^2 dxdt+\iint_{H_{j+1}^{1}}V(x)u^2(x,t)\psi^2dxdt\\
	&\leq \epsilon_2 \iint_{H_{j+1}^{1}} |\nabla u(x,t)|^2\psi^2 dxdt+\frac{1}{\epsilon_2}\iint_{H_{j+1}^{1}}| u(x,t)|^2|\nabla \psi|^2dxdt +\frac{1}{2}\iint_{H_{j+1}^{1}} u^2(x,t)\p_t \psi^2 dxdt.
	\end{aligned}
	$$
	By the assumption on $\psi$ and taking $\epsilon_2=\frac{1}{2}$, we can see,	
	\be\label{imp42}
	\iint_{H_{j}^{2}}|\nabla u(x,t)|^2dxdt+\iint_{H_{j}^{2}}V(x) u^2(x,t)dxdt \leq \frac{Ck}{t_0} \iint_{H_{j+1}^{1}}| u(x,t)|^2dxdt.
	\ee

	Combine the inequalities \eqref{imp41} and \eqref{imp42}, we have
	$$
	\iint_{H_j^1}|\p_tu(x,t)|^2dxdt\leq \frac{C^2k^2}{t_0^2} \iint_{H_{j+1}^1}| u(x,t)|^2dxdt.
	$$
	
	By Remark \ref{remark}, $\partial_t^l u $ is also a weak solution of \eqref{2ndeq} for any nonnegative integer $l$. Thence
	$$
	\iint_{H_{1}^{1}}\left(\partial_t^k u(x,t)\right)^2 d x d t\leq \frac{C^2k^2}{t_0^2}\iint_{H_{2}^{1}}\left(\partial_t^{k-1} u(x,t)\right)^2 d x d t\leq...\leq \frac{C^{2k} k^{2k}}{t_0^{2k}}\iint_{H_{k+1}^{1}} u(x,t)^2 d x d t.
	$$
	Therefore, by Lemma \ref{lem4new}, \eqref{volume1}
	$$
	\begin{aligned}
	|\partial_t^ku(x_0,t_0)|^2&\leq  C \left(\frac{k}{t_0}\right)^{d/2+1}e^{D d(x_0,0)}
	\iint_{Q_{k^{-\frac{1}{2}}}\left(x_{0},t_0\right)} |\partial_t^ku(x,t)|^2 d x d t\\
	&\leq  C \left(\frac{k}{t_0}\right)^{d/2+1} \left(\frac{Ck}{t_0}\right)^{2k}\iint_{H_{k+1}^1}\left(u(x,t)\right)^2dxdt\leq \frac{A_1^2A_5^{2k+2}k^{2k}}{t_0^{2k}}e^{4A_2d^2(x_0,0)}.
	\end{aligned}
	$$
	
	The rest of the proof is the same as Theorem \ref{2ndtheo}.
	\qed
	
	\begin{remark}
		To make sure the set of functions satisfying condition \ref{connew2} is nontrivial when $V(x)=\frac{A}{d(x,0)^2}$, we give some examples here. The first one is
		$$a_0(x)=\sum_{j=1}^{\infty}\frac{|x|^{2j}}{((2j)!)^{1+s}},$$ where $ s\geq 0$.
		Now we give a lemma explaining $a_0(x)$ satisfies condition \ref{connew2} in $R^d$. We can prove the following lemma by induction.
		
		\begin{lemma}
			Let the space $\mathrm{M}=R^d$, then there are two sequences of positive number $a_{j,k}$ and $b_{j,k}$ where $j,k$ are nonnegative integers satisfying $$\left(\Delta-\frac{2d}{|x|^2}\right)^{k}a_0(x)=\sum_{j=k+1}^{\infty}b_{j,k}|x|^{2j-2k}$$ and 
			$$\Delta^{k}a_0(x)=\sum_{j=k}^{\infty}a_{j,k}|x|^{2j-2k}.$$
			Besides, we have $0\leq b_{j,k}\leq a_{j,k},$ and 
			$$\bigg|\left(\Delta-\frac{2d}{|x|^2}\right)^{k}a_0(x)\bigg|\leq \Delta^{k}a_0(x) \leq C^k k! e^{4d |x|^2}.$$
		\end{lemma}
		
		Then we can have another example $a^*(x)=\sum_{j=1}^{\infty}\frac{(-1)^{j+1}|x|^{2j}}{((2j)!)^{1+s}}$, $s>0$ which also satisfies the condition \eqref{connew2}.
		This is because if we let
		$$\left(\Delta-\frac{2}{|x|^2}\right)^{k}a^*(x)=\sum_{j=k+1}^{\infty}d_{j,k}|x|^{2j-2k},$$
		then
		$d_{j,m+1}=(-1)^{j+1}b_{j,m+1}$ for any nonnegative integers $j$,$m$.
		
		Especially, we can also prove the functions $|x|^2cos(|x|)$ and $|x|sin(|x|)$ also satisfies the condition \eqref{connew2} by the same method.
	\end{remark}

	We have similar corollaries as Corollary \ref{cor1} and Corollary \ref{cor2} using the same proof.
	
	\begin{corollary}\label{cor3}
		Let $V=V(x)$ be a potential function satisfying either the conditions in Theorem \ref{4ththeo} or $V(x)\geq0$. Then
		the Cauchy problem for the backward heat equation with potentials
		$$
		\quad\left\{\begin{array}{l}
		\partial_{t} u(x,t)+(\Delta-V(x)) u(x,t)=0 \\
		u(x,0)=a(x),
		\end{array}\right.
		$$
		has a weak solution of exponential growth of order 2 in $\mathrm{M}\times (0,\delta)$ for some $\delta>0$ if and only if there exist some constants $A_{2}, A_3$ satisfying:
		$$
		\left|\left(\Delta-V(x)\right)^j a(x) \right| \leq A_{2}^{j+1} j^{j} e^{A_3 d^2(x,0)}, \quad j=0,1,2, \ldots
		$$
	\end{corollary}

	\begin{corollary}
		Let $V=V(x)$ satisfies the same conditions as Corollary \ref{cor3} above. Then the Cauchy problem
		$$
		\quad\left\{\begin{array}{l}
		\partial_{t} u(x,t)-(\Delta-V(x)) u(x,t)=0 \\
		u(x,0)=a(x)
		\end{array}\right.
		$$
		has a weak solution of exponential growth of order 2, which is also analytic in time in
		$\mathrm{M}\times[0, \delta) $ for some $\delta>0$ with a radius of convergence independent of $x$ if and only if there exist some constants $A_{2}, A_3$ satisfying:
		$$
		\left|\left(\Delta-V(x)\right)^j a(x) \right| \leq A_{2}^{j+1} j^{j} e^{A_3 d^2(x,0)}, \quad j=0,1,2, \ldots
		$$
	\end{corollary}

	\section{Nonlinear heat equations with power nonlinearity}\label{4thsec}
	This section is about some nonlinear heat equations with power nonlinearity of order $p$ \eqref{equa0} where $p\in(0,\infty)$. There are two main theorems \ref{theo00} and \ref{theo05} in this section and the main tools to prove them are Lemmas \ref{tri} and \ref{tri2}. We first prove the case when the solution $u$ is bounded and $p$ is an integer. Then we turn to the case when $0<C_3\leq|u|\leq C_4$ and $p$ is any rational number.

	For \eqref{equa0}, since we assume the solution $u$ is bounded, by standard theory, $u$ is actually smooth.
	We need a lemma about the time derivative of the heat kernel on $\mathrm{M}$ first.
	\begin{lemma}\label{heat1}
		Let $\mathrm{M}$ be the same manifold as Theorem \ref{theo00} above. Then for any $x,y\in \mathrm{M}$, $0<t\leq 1$ and any nonnegative integer $k$, there exist some constants $C_1$ and $C_5$ depending only on $\mathrm{M}$ and $d$ such that the heat kernel $\Gamma(x,t;y)$ of the heat equation
		$$\p_t u-\Delta u=0,$$
		satisfies the following condition:
		\be\label{heat}
		|\p_t^k \Gamma(x,t;y)|\leq \frac{C_1^{k+1}k^{k-2/3}}{t^{k}|B(x,\sqrt{t})|}e^{\frac{-C_5 d(x,y)^2}{t}}.
		\ee	
	\end{lemma}
	
	\begin{remark}
		To our best knowledge, up to now, in the literature, one just have
		$$
		|\p_t^k \Gamma(x,t;y)|\leq \frac{C(k)}{t^{k}|B(x,\sqrt{t})|}e^{\frac{-C_5 d(x,y)^2}{t}}
		$$ 
		in the manifold case, where $C(k)$ is not calculated explicitly. So here we obtain a more accurate result.
	\end{remark}
	\textbf{Proof of Lemma \ref{heat1}}\\[2mm]
	Fix any $t_0\in(0,1]$ and $x_0, y_0\in \mathrm{M}$, we would like to get the estimates of $\p_t^k \Gamma(x_0,t_0;y_0)$. For any nonnegative integer $k$ and $j=1,2,\cdots,k$,
	we define some space-time domains:
	$$M_j^1=\left\{(x,t):d(x,x_0)<\frac{j \sqrt{t_0}}{\sqrt{2k}},t\in \left(t_0-\frac{j t_0}{2k},t_0\right)\right\},$$
	$$M_j^2=\left\{(x,t):d(x,x_0)<\frac{(j+0.5) \sqrt{t_0}}{\sqrt{2k}},t\in \left(t_0-\frac{(j+0.5) t_0}{2k},t_0\right)\right\}.$$
	Then $M_j^1\subset M_j^2\subset M_{j+1}^1$.
	
	Following the method used in the proof of Theorem \ref{4ththeo}, for some constant $C$, it holds
	\be\label{im3}
	\iint_{M_1^1}|\p_t^k \Gamma(x,t;y_0)|^2 dxdt\leq \frac{C^{2k}k^{2k}}{t_0^{2k}}\iint_{M_{k+1}^1}|\Gamma(x,t;y_0)|^2 dxdt.
	\ee
	Then we need to use the well-known result for the upper bound of the heat kernel which can be found in P.Li$\&$S.T.Yau\cite{[LY]} or L.Saloff-Coste\cite{[Saloff]}, which is
	$$\Gamma(x,t;y)\leq \frac{C_3^\prime e^{\frac{-C_4^\prime d(x,y)^2}{t}}}{|B(x,\sqrt{t})|},\ \forall x,y\in \mathrm{M}\ \text{and}\ t\in(0,1],$$
	for some constants $C_3^\prime$ and $C_4^\prime$.

	Now we have two cases.\\
	\textbf{Case (1):} $\boldsymbol{d(y_0,x_0)\leq \sqrt{4kt_0}}$.\\ 
	In this case, using \eqref{volume1}
	$$
	\frac{C^{2k}k^{2k}}{t_0^{2k}}\iint_{M_{k+1}^1}|\Gamma(x,t;y_0)|^2 dxdt\leq \frac{C^{2k+1/2}k^{2k}e^{\frac{D(k+1)^2t_0}{2k t_0}}}{t_0^{2k-1}|B(x_0,\sqrt{t_0})|}\leq \frac{C^{2k+1}k^{2k+1}}{t_0^{2k-1}|B(x_0,\sqrt{t_0})|}e^{\frac{-C_5 d(x_0,y_0)^2}{t_0}},
	$$
	for some constant $C$.
	
	\textbf{Case (2):} $\boldsymbol{d(y_0,x_0)>\sqrt{4kt_0}}$.\\ 
	In this case, because $d(x,x_0)<\frac{(k+1)\sqrt{t_0}}{\sqrt{2k}}$, $\frac{\sqrt{2}-1}{\sqrt{2}}< \frac{d(x,y_0)}{d(x_0,y_0)}< 2$. Therefore,
	$$
	\begin{aligned}
	&\frac{C^{2k}k^{2k}}{t_0^{2k}}\iint_{M_{k+1}^1}|\Gamma(x,t;y_0)|^2 dxdt\\
	&\leq \frac{C^{2k}k^{2k}t_0|B(x_0,\frac{(k+1)\sqrt{t_0}}{\sqrt{2k}})|e^{\frac{2D(k+1)^2 t_0}{ 2k t_0}}}{t_0^{2k}|B(x_0,\sqrt{t_0})|^2}e^{\frac{-(3-2\sqrt{2})C_4^\prime d(x_0,y_0)^2}{2t_0}}\\
	&\leq \frac{C^{2k+1/2}k^{2k+1}}{t_0^{2k-1}|B(x_0,\sqrt{t_0})|}e^{\frac{-C_5 d(x_0,y_0)^2}{t_0}}\leq \frac{C^{2k+1}k^{2k+1}}{t_0^{2k-1}|B(x_0,\sqrt{t_0})|}e^{\frac{-C_5 d(x_0,y_0)^2}{t_0}}.
	\end{aligned}
	$$
	
	Combine the above two cases,
	\be\label{im1}
	\iint_{M_1^1}|\p_t^k \Gamma(x,t;y_0)|^2 dxdt\leq \frac{C^{2k+1}k^{2k+1}}{t_0^{2k-1}|B(x_0,\sqrt{t_0})|}e^{\frac{-C_5 d(x_0,y)^2}{t_0}}.
	\ee
	Then we recall a well-known parabolic mean value inequality which can be found, for instance, in Theorem 14.7 of P.Li\cite{[Li]}. To be more precise, by the method of Lemma \ref{lem3new}, for any $0<p<\infty$ and $0<r<R<1$, any nonnegative subsolution $u=u(x,t)$ of the heat equation satisfies
	$$\begin{aligned}
	&\sup _{Q_{r}\left(x_{0},t_0 \right)} u(x,t)^p\leq C\left(\frac{R^2}{|B(x_0,R)|^{\frac{2}{\nu_2}}}\right)^{\frac{1}{\theta^*-1}}
	\left(  \frac{1}{|R-r|^2}\right)^{\frac{\theta^*}{\theta^*-1}}\iint_{Q_{R}\left(x_{0},t_0\right)} u(x,t)^p d x d t ,
	\end{aligned}$$
	where $\theta^*=1+\frac{2}{\nu_2}$ and $\nu_2$ is defined in \eqref{sobolevfinal1}. Let $u(x,t)=|\p_t^k \Gamma(x,t;y_0)|^2$, $p=1$, $r=0$ and $R=\sqrt{t_0} / \sqrt{2k}$, we can see
	\be\label{meanva}
	\begin{aligned}
		|\p_t^k \Gamma(x_0,t;y_0)|^{2}&\leq \frac{C k}{\left|B\left(x_{0}, \sqrt{t_0} / \sqrt{2k}\right)\right|t_0} \iint_{Q_{\sqrt{t_0} / \sqrt{2k}}\left(x_{0},t_0\right)} (\p_t^k \Gamma(x,t;y_0))^{2}d x d t\\
		& \leq \frac{C k^{d/2+1}}{\left|B\left(x_{0}, \sqrt{t_0}\right)\right|t_0} \iint_{Q_{\sqrt{t_0} / \sqrt{2k}}\left(x_{0},t_0\right)} (\p_t^k \Gamma(x,t;y_0))^{2}d x d t,
	\end{aligned}
	\ee
	where we have used the Bishop-Gromov volume comparison theorem in the last inequality. \\
	By \eqref{im3},\eqref{im1} and \eqref{meanva}, we see 
	$$
	(\p_t^k \Gamma(x_0,t_0;y_0))^{2}\leq \frac{C^{2k+2}k^{2k+d/2+2}}{t_0^{2k}|B(x_0,\sqrt{t_0})|^2}e^{\frac{-C_5 d(x_0,y_0)^2}{t_0}}.
	$$
	Thus, $$|\p_t^k \Gamma(x_0,t_0;y_0)|\leq \frac{C_1^{k+1}k^{k-2/3}}{t_0^{k}|B(x_0,\sqrt{t_0})|}e^{\frac{-C_5 d(x_0,y_0)^2}{t_0}},$$
	for some $C_1$ large enough, which finishes the proof of Lemma \ref{heat1}.
	\qed
	
	\begin{remark}
		By the estimate of the time derivative of heat kernel $\Gamma(x,t;y)$, we can see the solution $u=u(x,t)$ of heat equation $u_t-\Delta u=0$ is analytic in time if $u$ is of exponential growth of order 2 directly.
	\end{remark}

	Let $\binom{n}{i_1,i_2,\cdots, i_k}:=\frac{n!}{i_1!i_2!\cdots(n-i_1-i_2-\cdots-i_k)!}$. Then we have a lemma which will be used frequently.
	\begin{lemma}\label{tri}
		For any integers $n>1$ and $k>1$, there exists some constant $C=C(k)$ such that,
		$$
		\begin{aligned}
		&\sum_{\Sigma_{m=1}^k i_m<n,i_m>0}\binom{n}{i_1,i_2,\cdots, i_k}i_1^{i_1-2/3}i_2^{i_2-2/3}\cdots(n-i_1-i_2-\cdots-i_k)^{n-i_1-i_2-\cdots-i_k-2/3}\\
		&\leq Cn^{n-2/3}.
		\end{aligned}
		$$
	\end{lemma}
	This lemma is just an extension of the Lemma 3.2 of H.Dong$\&$ Q.Zhang\cite{[Zhang]} and we can prove it by the induction method and the Stirling formula.\\
	\pf
	$$
	\begin{aligned}
	&\sum_{\Sigma_{m=1}^k i_m<n,i_m>0}\binom{n}{i_1,i_2,\cdots, i_k}i_1^{i_1-2/3}i_2^{i_2-2/3}\cdots(n-i_1-i_2-\cdots-i_k)^{n-i_1-i_2-\cdots-i_k-2/3}\\
	&=\sum_{i_1=1}^{n}\binom{n}{i_1}i_1^{i_1-2/3}\sum_{\Sigma_{m=2}^k i_m<n-i_1,i_m>0}\binom{n-i_1}{i_2,\cdots, i_k}i_2^{i_2-2/3}\cdots(n-i_1-i_2-\cdots-i_k)^{n-i_1-i_2-\cdots-i_k-2/3}\\
	&\leq C \sum_{i_1=1}^{n-1}\binom{n}{i_1}i_1^{i_1-2/3}(n-i_1)^{n-i_1-2/3}\leq Cn^{n-2/3}\sum_{i_1=1}^{n-1}\frac{n^{7/6}}{i_1^{7/6}(n-i_1)^{7/6}}\\
	&\leq Cn^{n-2/3}\sum_{i_1=1}^{n-1}\left(\frac{1}{i_1}+\frac{1}{n-i_1}\right)^{7/6}\leq Cn^{n-2/3}.
	\end{aligned}
	$$
	\qed\\

	Then we have the following lemma to connect $\p_t^n(t^n u^p)$ and $\p_t^n(t^n u)$ for any positive integer $n$.
	\begin{lemma}\label{tri2}
		Let $f_1(t)$,$f_2(t)$,$\cdots$,$f_k(t)$ be smooth functions. For any nonnegative integer $n$, we have
		$$
		\begin{aligned}
		&\p_t^n(t^n f_1(t)f_2(t)\cdots f_k(t))\\
		=&\sum\limits_{m=0}^{k-1}(-1)^m \frac{n!}{(n-m)!}\binom{k-1}{m}\sum\limits_{i_l\geq 0}\binom{n-m}{i_1,i_2,\cdots, i_{k-1}}\\
		&\p_t^{i_1}(t^{i_1}f_1(t))\cdots \p_t^{i_{k-1}}(t^{i_{k-1}}f_{k-1}(t))\p_t^{n-m-\Sigma_{l=1}^{k-1}i_l}(t^{n-m-\Sigma_{l=1}^{k-1}i_l}f_{k}(t)).
		\end{aligned}
		$$
		Here for $\binom{n}{i_1,i_2,\cdots,i_k}$ we always assume $\sum\limits_{l=1}^k i_l\leq n$. 
	\end{lemma}

	\pf We can prove it by induction using Lemma 3.3 of H.Dong$\&$Q.Zhang\cite{[Zhang]}.

	\begin{remark}
		Especially, when $f_1=f_2=\cdots =f_k=f$, it holds
		\be\label{1stpart}
		\begin{aligned}
			&\p_t^n(t^n f^k(t))\\
			=&\sum\limits_{m=0}^{k-1}(-1)^m \frac{n!}{(n-m)!}\binom{k-1}{m}\sum\limits_{i_l\geq 0}\binom{n-m}{i_1,i_2,\cdots, i_{k-1}}\\
			&\p_t^{i_1}(t^{i_1}f(t)) \cdots \p_t^{i_{k-1}}(t^{i_{k-1}}f(t))\p_t^{n-m-\Sigma_{l=1}^{k-1}i_l}(t^{n-m-\Sigma_{l=1}^{k-1}i_l}f(t)).
		\end{aligned}
		\ee
		Moreover, when $f_i(t)=f(t)^{\frac{1}{k}}$ for any $i=1,\cdots, k$, we have
		\be\label{2ndpart}
		\begin{aligned}
			&k f(t)^{\frac{k-1}{k}}\p_t^n(t^n f(t)^{\frac{1}{k}})\\
			&=\p_t^n(t^n f(t))-\sum\limits_{m=1}^{k-1}(-1)^m \frac{n!}{(n-m)!}\binom{k-1}{m}\sum\limits_{i_l\geq 0}\binom{n-m}{i_1,i_2,\cdots, i_{k-1}}\\
			&\p_t^{i_1}(t^{i_1}f(t)^{\frac{1}{k}}) \cdots \p_t^{i_{k-1}}(t^{i_{k-1}}f(t)^{\frac{1}{k}})\p_t^{n-m-\Sigma_{l=1}^{k-1}i_l}(t^{n-m-\Sigma_{l=1}^{k-1}i_l}f(t)^{\frac{1}{k}})\\
			&\quad-\sum\limits_{\mbox{\tiny$\begin{array}{c}
					n>i_l\geq 0\\
					\Sigma_{l=1}^{k-1}i_l>0\end{array}$}}\binom{n}{i_1,i_2,\cdots, i_{k-1}}\p_t^{i_1}(t^{i_1}f(t)^{\frac{1}{k}}) \cdots \p_t^{i_{k-1}}(t^{i_{k-1}}f(t)^{\frac{1}{k}})\p_t^{n-\Sigma_{l=1}^{k-1}i_l}(t^{n-\Sigma_{l=1}^{k-1}i_l}f(t)^{\frac{1}{k}}).
		\end{aligned}
		\ee
	\end{remark}

	We first establish the following proposition before embarking on the proof of Theorem \ref{theo00}.

	\begin{proposition}
		Under the conditions of Theorem \ref{theo00} above, for any integer $n \geq 1$, it holds
		\be\label{fin} \left\|\partial_{t}^{n}\left(t^{n} u(\cdot,t)\right)\right\|_{L^{\infty}\left(\mathrm{M}\right)} \leq N^{n-1 / 2} n^{n-2 / 3} 
		\ee
		for some sufficiently large constant  $N \geq 1$. \\ 
	\end{proposition}
	\pf By induction and by lemma \ref{heat}, there exist some constant $C_1$ such that for any integer $k>1$,
	$$\left\|\partial_{t}^{k}\left(t^{k} \Gamma( \cdot,t)\right)\right\|_{L^1(\mathrm{M})} \leq C_{1}^{k+1} k^{k-2 / 3}.$$
	We shall prove the proposition inductively. As $u$  is a solution, we have 
	$$u(x,t)=\int_{\mathrm{M}}\Gamma(x,t;y) u(0, y)dy+\int_{0}^{t} \int_{\mathrm{M}}\Gamma(x,t-s;y) u^p (y,s)dy d s,$$
	as a consequence,
	\be\label{u}
	\begin{aligned} 
		\partial_t^n(t^{n}u(x,t))=& \int_{\mathrm{M}}\partial_{t}^{n}(t^{n}\Gamma(x,t;y))  u(0, y)dy+\partial_{t}^{n}(\int_{\mathrm{M}}\int_{0}^{t}t^{n} \Gamma(x,t-s;y)  u^p(y,s) dyds) \\&
		:=I_{1}+I_{2}. 
	\end{aligned}
	\ee
	It holds
	\be\label{I_1fi}|I_{1}| \leq C_{2} C_{1}^{n+1} n^{n-2 / 3} \leq N^{n-2 / 3} n^{n-2 / 3}\ee
	for sufficiently large $N $.\\
	
	To estimate $I_2$, similar to the inequality (3.7) from the paper H.Dong$\&$ Q.Zhang\cite{[Zhang]}, we yield
	\be\label{I_21}
	\begin{aligned}
		I_{2} &=\sum_{k=0}^{n} \binom{n}{k}\partial_{t}^{n} \int_{0}^{t}\int_{\mathrm{M}}\left((t-s)^{k} \Gamma(x,t-s;y)\right) \left(s^{n-k}u^p(y,s)\right) dyd s \\
		&=\sum_{k=0}^{n}\binom{n}{k} \partial_{t}^{n-k} \int_{0}^{t} \int_{\mathrm{M}} \partial_{t}^{k}\left((t-s)^{k} \Gamma(x,t-s;y)\right) \left(s^{n-k}u^p(y,s)\right) dyd s \\
		&=\sum_{k=0}^{n}\binom{n}{k} \partial_{t}^{n-k}\int_{0}^{t} \int_{\mathrm{M}}\partial_{s}^{k}\left(s^{k} \Gamma(x,s;y)\right) \left((t-s)^{n-k}u^p(y,t-s)\right) dyd s \\
		&=\sum_{k=0}^{n}\binom{n}{k} \int_{0}^{t} \int_{\mathrm{M}}\partial_{s}^{k}\left(s^{k} \Gamma(x,s;y)\right)  \partial_{t}^{n-k}\left((t-s)^{n-k}u^p(y,t-s)\right) dyd s.
	\end{aligned}
	\ee
	
	Using Lemma \ref{tri} and equality \ref{1stpart}, the it holds by induction
	$$|\partial_{t}^{n}\left(t^{n}(u^p(x,t))\right)| \leq pC_2^{p-1}\left|\partial_{t}^{n}\left(t^{n} u(x,t)\right)\right|+N^{n-3 / 4} n^{n-2 / 3},$$
	and for $k=1, \ldots, n-1$
	$${\left|\partial_{t}^{k}\left(t^{k}(u^p(x,t))\right)\right|}  { \leq N^{k-1 / 3} k^{k-2 / 3}}.$$
	Following the similar procedure as in the paper H.Dong$\&$ Q.Zhang\cite{[Zhang]}, we have 
	\be\label{I_2fi}
	\begin{aligned}
		\left|I_{2}\right| \leq & \int_{0}^{t}C_{1}^{n+1} n^{n-2/3} C_{2}^{p}+C\left(pC_{2}^{p-1}\left\|\partial_{t}^{n}\left((t-s)^{n} u( \cdot,t-s)\right)\right\|_{L^{\infty}}+N^{n-3 / 4} n^{n-2 / 3}\right)\\
		&+\sum_{k=1}^{n-1}\binom{n}{k} C_{1}^{k+1} k^{k-2 / 3} \cdot N^{n-k-1 / 3}(n-k)^{n-k-2 / 3} d s \\
		\leq & N^{n-2 / 3} n^{n-2 / 3}t+ CpC_2^{p-1} \int_{0}^{t} \left\|\partial_{s}^{n}\left(s^{n} u(\cdot,s)\right)\right\|_{L^{\infty}} d s
	\end{aligned}
	\ee
	for sufficiently large $N$ depending on $C_{1}$, $C_{2}$,$p$, $d$ and $K_0$.
	
	Combining the estimates of $I_1$ \eqref{I_1fi} and $I_2$ \eqref{I_2fi}, we can get \eqref{fin} by applying Gronwall's inequality and finish the proof of the proposition.
	\qed\\
	Now we begin the proof of the theorem \ref{theo00}. \\
	\subsection{Proof of Theorem \ref{theo00}}
	This part is the same as H.Dong$\&$ Q.Zhang\cite{[Zhang]}. We just copy it down here for the convenience of reading.
	
	Note that
	$$
	\partial_{t}^{n}\left(t^{k} u\right)=n \partial_{t}^{n-1}\left(t^{k-1} u\right)+t \partial_{t}^{n}\left(t^{k-1} u\right).
	$$
	Taking $k=n$, we obtain
	$$
	\sup _{t \in(0,1]}\left\|t \partial_{t}^{n}\left(t^{n-1} u( \cdot,t)\right)\right\|_{L^{\infty}\left(\mathrm{M}\right)} \leq N^{n}(1+1 / N) n^{n}.
	$$
	By induction,
	$$
	\sup _{t \in(0,1]}\left\|t^{n} \partial_{t}^{n} u(\cdot,t)\right\|_{L^{\infty}\left(\mathrm{M}\right)} \leq N^{n}(1+1 / N)^{n} n^{n}=(N+1)^{n} n^{n}.
	$$
	The theorem is proved.
	\qed\\
	
	To prove Theorem \ref{theo05}, we also have a proposition first using Lemmas \ref{tri} and \ref{tri2}.

	\begin{proposition}\label{prp}
		Under the conditions of Theorem \ref{theo05} above, for any integer $n \geq 1$, we have
		$$ \big|\partial_{t}^{n}\left(t^{n} u(x,t)\right)\big|\leq N^{n-1 / 2} n^{n-2 / 3},$$
		for some sufficiently large constant $N$.
	\end{proposition}
	\pf We shall prove the proposition inductively. First, we can get equality \eqref{u} in the same way. Then similar to inequality \eqref{I_1fi}, we see
	$$
	|I_{1}| \leq  N^{n-2 / 3} n^{n-2 / 3},
	$$
	for sufficiently large $N $.
	
	By equality \eqref{2ndpart} and Lemma \ref{tri}, we can prove
	by induction, for any $k=1,2,\cdots, n-1$\\
	$$\big|\p_t^k(t^ku(x,t)^{1/q_2})\big|\leq N^{k-5/12} k^{k-2 / 3},$$
	and 
	$$\bigg|\frac{\p_t^n(t^nu(x,t)^{1/q_2})}{u(x,t)^{1/q_2}}\bigg|\leq \frac{1}{q_2}\bigg|\frac{\p_t^n(t^n u(x,t))}{ u(x,t)}\bigg|+ N^{n-19/24}n^{n-2/3}.$$
	To be more precise, if we assume for any $l=1,2,\cdots, k-1$\\
	$$\big|\p_t^l(t^lu(x,t)^{1/q_2})\big|\leq N^{l-5/12} l^{l-2 / 3},$$ then
	$$
	\begin{aligned}
	&k C_3^{\frac{q_2-1}{q_2}}|\p_t^k(t^k u(x,t)^{\frac{1}{q_2}})|\\
	&\leq |\p_t^k(t^k u(x,t))|+\sum\limits_{m=1}^{q_2-1} \frac{k!}{(k-m)!}\binom{q_2-1}{m}\sum\limits_{i_l\geq 0}\binom{k-m}{i_1,i_2,\cdots, i_{q_2-1}}\\
	&N^{i_1-5/12} i_1^{i_1-2 / 3} \cdots N^{i_{q_2-1}-5/12} i_{q_2-1}^{i_{q_2-1}-2 / 3}N^{k-m-\Sigma_{l=1}^{q_2-1}i_l-5/12}(k-m-\Sigma_{l=1}^{q_2-1}i_l)^{k-m-\Sigma_{l=1}^{q_2-1}i_l}\\
	&\quad+\sum\limits_{\mbox{\tiny$\begin{array}{c}
			k>i_l\geq 0\\
			\Sigma_{l=1}^{q_2-1}i_l>0\end{array}$}}\binom{k}{i_1,i_2,\cdots, i_{q_2-1}}N^{i_1-5/12} i_1^{i_1-2 / 3} \cdots  N^{i_{q_2-1}-5/12} i_{q_2-1}^{i_{q_2-1}-2 / 3}\\
	&N^{k-\Sigma_{l=1}^{q_2-1}i_l-5/12}(k-m-\Sigma_{l=1}^{q_2-1}i_l)^{k-m-\Sigma_{l=1}^{q_2-1}i_l}\\
	&\leq |\p_t^m(t^m u(x,t))|+N^{k-1/2}.
	\end{aligned}
	$$
	
	Therefore by equality \eqref{1stpart} and Lemma \ref{tri}, we can prove
	by induction that for any  $k=1,2,\cdots, n-1$\\
	$$\big|\p_t^k(t^ku(x,t)^{q_1/q_2})\big|\leq N^{k-1/3} k^{k-2 / 3},$$
	and 
	$$\bigg|\frac{\p_t^n(t^nu(x,t)^{q_1/q_2})}{u(x,t)^{q_1/q_2}}\bigg|\leq q_1\bigg|\frac{\p_t^n(t^n u(x,t)^{1/q_2})}{u(x,t)^{1/q_2}}\bigg|+ N^{n-3/4}n^{n-2/3}\leq \frac{q_1}{q_2}\bigg|\frac{\p_t^n(t^n u(x,t))}{ u(x,t)}\bigg|+ N^{n-3/4}n^{n-2/3},$$
	for some constant $N$ large enough.\\
	
	Therefore by \eqref{I_21},
	$$
	\begin{aligned}
	\left|I_{2}\right| \leq & \int_{0}^{t}C_4 C_{1}^{n+1} n^{n-2/3}+ C\left(\frac{q_1}{q_2} \bigg|\frac{\p_t^n((t-s)^n u(\cdot,t-s))}{ u(\cdot,t)^{1-q_1/q_2}}\bigg|_{L^\infty(\mathrm{M})}+ C_4^{q_1/q_2}N^{n-3/4}n^{n-2/3}\right)\\
	&+\sum_{k=1}^{n-1}\binom{n}{k} C_{1}^{k+1} k^{k-2 / 3} \cdot N^{n-k-1 / 3}(n-k)^{n-k-2 / 3} d s \\
	\leq & N^{n-2 / 3} n^{n-2 / 3}t+ C \int_{0}^{t} \left\|\partial_{s}^{n}\left(s^{n} u(\cdot,s)\right)\right\|_{L^{\infty}} ds,
	\end{aligned}
	$$
	for sufficiently large $N$ depending on $C_1$, $C_{3}$, $C_{4}$, $p$, $d$ and $K_0$. Using the estimates of $I_1$, $I_2$ above and Gronwall's inequality, we can finish the proof of Proposition \ref{prp}.
	\qed
	
	With this proposition at hand, we can prove the Theorem \ref{theo05} immediately.
	
	\subsection{Proof of Theorem \ref{theo05} }
	The proof is exactly the same as the last part of the proof of Theorem \ref{theo00}.
	\qed
	
	\begin{remark}
		For the case when $0<p<1$, we can have a particular solution 
		$$
		\quad
		u(x,t)=\left\{\begin{array}{l}
		\left((1-p)(t-\frac{1}{2})\right)^{\frac{1}{1-p}} \quad \text{when} \quad \frac{1}{2}<t<1\\
		0 \quad\quad\text{when} \quad0\leq t\leq\frac{1}{2},
		\end{array}\right.
		$$
		which is not analytic at $t=\frac{1}{2}$. We can use this example to say that $u$ may not be allowed to be $0$ to get the time analyticity conclusion.
	\end{remark}
	\begin{remark}
		For the time analyticity at $t=0$, according to the paper 
		G.Lysik$\&$S.Michalik\\\cite{[Lysik]}, even for some polynomial functions $f(u)$, the formal solutions for  $\p_tu(x,t)-\Delta u(x,t)=f(u)$ are not in general analytic at $t=0$ even if the initial condition is analytic.
	\end{remark}	
	\begin{remark}
		It is maybe true that the conclusion in Theorem \ref{theo05} can be extended to all the real number $p$.
	\end{remark}
	\section*{Acknowledgement}
	The author wishes to express his appreciation to his advisor Professor Qi S. Zhang for providing him with this
	problem, sharing 
	the ideas and offering a lot of helpful discussions and suggestions. Besides, the author feels grateful to Professor Hongjie Dong, Professor Xin Yang
	,Professor Na Zhao for going over the paper and making suggestions.


\begin{thebibliography}{99}
		\bibitem{[BD]} Barbatis, Gerassimos, Davies, E.. (1997). Sharp Bounds on Heat Kernels of Higher Order Uniformly Elliptic Operators. J. Operator Theory. 36. 
		\bibitem{[CC]}	Cheeger, J. and Colding, T. H. (1996). Lower bounds on Ricci curvature and the almost rigidity of warped products. Annals of Mathematics, 144(1), 189-237. 
		\bibitem{[DK]}Dong, Hongjie; Kim, Doyoon On the Lp-solvability of higher order parabolic and elliptic systems with BMO coefficients. Arch. Ration. Mech. Anal. 199 (2011), no. 3, 889–941.
		\bibitem{[DP]} Dong, Hongjie, Pan, Xinghong.
		Time analyticity for inhomogeneous parabolic equations and the Navier-Stokes equations in the half space. 
		J. Math. Fluid Mech. 22 (2020).
		\bibitem{[Zhang]} Dong, Hongjie and Zhang, Qi. (2020). Time analyticity for the heat equation and Navier-Stokes equations. Journal of Functional Analysis.
		\bibitem{[LE]}L. Escauriaza, S. Montaner and C. Zhang, Analyticity of solutions to parabolic evolutions
		and applications, SIAM J. Math. Anal. 49 (2017), no. 5, 4064–4092.
		\bibitem{[Giaquinta]}Giaquinta, M., Introduction to Regularity Theory for Nonlinear Elliptic
		Systems, Basel, Boston, Berlin: Birkh$\ddot{a}$user, 1993.
		\bibitem{[GI]} Giga, Yoshikazu, Time and spatial analyticity of solutions of the Navier-Stokes equations, Comm.
		Partial Differential Equations 8 (1983), no. 8, 929-948.
		\bibitem{[Hebey]}Emmanuel Hebey. Nonlinear Analysis on Manifolds: Sobolev Spaces and Inequalities,
		Courant Lecture Notes in Mathematics, Vol. 5 (New York University, Courant Institute
		of Mathematical Sciences, New York; American Mathematical Society, Providence,
		RI, 1999)
		\bibitem{[HL]}Qi Hou, Laurent Saloff-Coste. Time regularity for local weak solutions of the heat equation on local Dirichlet spaces, 	arXiv:1912.12998
		\bibitem{[Komatsu]}G. Komatsu. Global analyticity up to the boundary of solutions of the Navier-Stokes equation,
		Comm. Pure Appl. Math. 33 (1980), no. 4, 545–566.
		\bibitem{[Saloff]}Saloff-Coste,L. (2001). Aspects of Sobolev-Type Inequalities (London Mathematical Society Lecture Note Series). Cambridge: Cambridge University Press. doi:10.1017/CBO9780511549762
		\bibitem{[Li]} Li, Peter, Geometric analysis. Cambridge Studies in Advanced Mathematics, 134. Cambridge Univer-
		sity Press, Cambridge, 2012. x+406 pp.
		\bibitem{[LP]} Li, Z., Pan, X.: Some remarks on regularity criteria of axially symmetric Navier–Stokes equations. Commun. Pure Appl. Anal. 18(3), 1333–1350 (2019)
		\bibitem{[LZ]}F. Lin and Q. S. Zhang, On ancient solutions of the heat equations, Comm. Pure Appl. Math.
		72 (2019), no. 9, 2006–2028.
		\bibitem{[ZZ]}Zijin Li, Qi S. Zhang.
		Regularity of weak solutions of elliptic and parabolic equations with some critical or supercritical potentials,
		Journal of Differential Equations,
		Volume 263, Issue 1,
		2017,
		Pages 57-87.
		\bibitem{[LY]}Li, Peter; Yau, Shing-Tung
		On the parabolic kernel of the Schrödinger operator.
		Acta Math. 156 (1986), no. 3-4, 153–201.
		58G11 (35J10)
		\bibitem{[Lysik]}Grzegorz \L ysik, S \L awomir Michalik. Formal solutions of semilinear heat equations, Journal of Mathematical Analysis and Applications, 2007. 
		\bibitem{[Masuda]}K. Masuda, On the analyticity and the unique continuation theorem for solutions of the
		Navier-Stokes equation, Proc. Japan Acad. 43 (1967), 827–832.
		\bibitem{[B]} Buser, Peter. A note on the isoperimetric constant. Annales scientifiques de l$'$\'Ecole Normale 
		Sup\'erieure, Serie 4, Volume 15 (1982) $no. 2, pp. 213-230. doi: 10.24033/asens.1426.$
		\bibitem{[Pr]} Promislow, Keith, Time analyticity and Gevrey regularity for solutions of a class of dissipative partial
		differential equations, Nonlinear Anal. 16 (1991), no. 11, 959-980.
		\bibitem{[Vas]} Vassilis G. Papanicolaou, Eva Kallitsi, George Smyrlis, Analytic Solutions of the Heat Equation, 2019.
		\bibitem{[SY]}Schoen, Richard (Richard M.), and Shing-Tung Yau. Lectures on Differential Geometry . Cambridge, MA: International Press, 1994. Print.
		\bibitem{[Tao]} Tao, Terence and Rodgers, Brad. (2018). The De Bruijn-Newman constant is non-negative. Forum of Mathematics, Pi. 8. 10.1017/fmp.2020.6. 
		\bibitem{[Widder]} D. V. Widder. Analytic solutions of the heat equation, Duke Math. J. 29 (1962), 497–503.
		\bibitem{[WZ]}B. Wong and Qi S. Zhang. Refined gradient bounds, possion equations and some applications
		to open k$\ddot{a}$hler manifolds. Asian J. Math, 7(3):1–28, September 2003.
		\bibitem{[Zhang2]} Qi Zhang.  A note on time analyticity for ancient solutions of the heat equation,  2019. Proc. Amer. Math. Soc. 148 (2020), 1665-1670
		\bibitem{[Zhangbook]} Zhang, Q.S.: Sobolev inequalities, heat kernels under Ricci flow, and the Poincar\'e conjecture. CRC Press,
		Boca Raton (2011).
		
	\end{thebibliography}
\end{document}